\documentclass[a4paper,11pt]{amsart}
\usepackage{amssymb,amsmath}
\usepackage{amsopn}
\usepackage{mathrsfs}
\usepackage{color}
\usepackage{wasysym}
\usepackage{vmargin}
\usepackage{vmargin}
\usepackage{color}
\usepackage{amscd}
\usepackage{verbatim}
\newtheorem{theorem}{Theorem}[section]
\newtheorem{lemma}[theorem]{Lemma}

\numberwithin{equation}{section}

\newcommand\RR{{{\mathbb R}}}

\newcommand\SSS{{\mathbb S}}
\newcommand{\rr}{\mathbb{R}}
\newcommand{\eps}{\varepsilon}
\newcommand{\nn}{\mathbb{N}}

\def\un{{\mathrm{1~\hspace{-1.4ex}l}}}

\def\Id{\operatorname{Id}}

\def\N{\mathbb N}

\def\R{\mathbb R}

\def\val#1{\vert#1\vert}

\def\l2{L^2(\R^{n})}
\def\L2{L^2(\R^{2n})}

\def\vs{\vskip.3cm}

\def\mat22#1#2#3#4{\begin{pmatrix}#1&#2\\ #3&#4\end{pmatrix}}

\begin{document}

\title[Gelfand-Shilov smoothing effect for the Boltzmann equation]
{Gelfand-Shilov smoothing properties of the radially symmetric spatially homogeneous Boltzmann equation without angular cutoff}
\author{N. Lerner, Y. Morimoto, K. Pravda-Starov  \& C.-J. Xu}
\date{\today}
\address{\noindent \textsc{N. Lerner, Institut de Math\'ematiques de Jussieu,
Universit\'e Pierre et Marie Curie (Paris VI),
4 Place Jussieu,
75252 Paris cedex 05,
France}}
\email{lerner@math.jussieu.fr}
\address{\noindent \textsc{Y. Morimoto, Graduate School of Human and Environmental Studies,
Kyoto University, Kyoto 606-8501, Japan}}
\email{morimoto@math.h.kyoto-u.ac.jp }
\address{\noindent \textsc{K. Pravda-Starov,
Universit\'e de Cergy-Pontoise,
D\'epartement de Math\'ematiques, CNRS UMR 8088,
95000 Cergy-Pontoise, France}}
\email{karel.pravda-starov@u-cergy.fr}
\address{\noindent \textsc{C.-J. Xu, School of Mathematics, Wuhan university 430072, Wuhan, P.R. China\\
 and  \\
 Universit\'e de Rouen, CNRS UMR 6085, D\'epartement de Math\'ematiques, 76801 Saint-Etienne du Rouvray, France}}
\email{Chao-Jiang.Xu@univ-rouen.fr}

\begin{abstract}
We prove that the Cauchy problem associated to the radially symmetric spatially homogeneous non-cutoff Boltzmann equation with Maxwellian molecules enjoys the same Gelfand-Shilov regularizing effect as the Cauchy problem defined by the evolution equation associated to a fractional harmonic oscillator.
\end{abstract}

\maketitle

\section{Introduction}
\subsection{Gelfand-Shilov regularity}
The Cauchy problem defined by the evolution equation associated to the harmonic oscillator
\begin{equation}\label{gel1}
\begin{cases}
\partial_tf+\mathcal{H}f=0,\\
f|_{t=0}=f_0 \in L^2(\rr^d),
\end{cases}\qquad \mathcal{H}=-\Delta_v+\frac{|v|^2}{4},
\end{equation}
enjoys nice regularizing properties. The smoothing effect for the solutions to this Cauchy problem is naturally described in term of the Gelfand-Shilov regularity~ \cite{gelfand,rodino1,rodino,toft}.
The Gelfand-Shilov spaces $S_{\nu}^{\mu}(\rr^d)$, with $\mu,\nu>0$, $\mu+\nu\geq 1$, are defined as the spaces of smooth functions $f \in C^{\infty}(\rr^d)$ satisfying the estimates
$$\exists A,C>0, \quad |\partial_v^{\alpha}f(v)| \leq C A^{|\alpha|}(\alpha !)^{\mu}e^{-\frac{1}{A}|v|^{1/\nu}}, \quad v \in \rr^d, \ \alpha \in \nn^d,$$
or, equivalently
$$\exists A,C>0, \quad \sup_{v \in \rr^d}|v^{\beta}\partial_v^{\alpha}f(v)| \leq C A^{|\alpha|+|\beta|}(\alpha !)^{\mu}(\beta !)^{\nu}, \quad \alpha, \beta \in \nn^d.$$
These Gelfand-Shilov spaces  $S_{\nu}^{\mu}(\rr^d)$ may also be characterized as the spaces of Schwartz functions $f \in \mathscr{S}(\rr^d)$ satisfying the estimates
$$\exists C>0, \eps>0, \quad |f(v)| \leq C e^{-\eps|v|^{1/\nu}}, \quad v \in \rr^d, \qquad |\widehat{f}(\xi)| \leq C e^{-\eps|\xi|^{1/\mu}}, \quad \xi \in \rr^d.$$
In particular, we notice that Hermite functions belong to the symmetric Gelfand-Shilov space  $S_{1/2}^{1/2}(\rr^d)$. More generally, the symmetric Gelfand-Shilov spaces $S_{\mu}^{\mu}(\rr^d)$, with $\mu \geq 1/2$, can be nicely characterized through the decomposition into the Hermite basis $(\Psi_{\alpha})_{\alpha \in \nn^d}$, see e.g. \cite{toft} (Proposition~1.2),
\begin{multline}\label{gel4}
f \in S_{\mu}^{\mu}(\rr^d) \Leftrightarrow f \in L^2(\rr^d), \ \exists t_0>0, \ \big\|\big((f,\Psi_{\alpha})_{L^2}\exp({t_0|\alpha|^{\frac{1}{2\mu}})}\big)_{\alpha \in \nn^d}\big\|_{l^2(\nn^d)}<+\infty\\
\Leftrightarrow f \in L^2(\rr^d), \ \exists t_0>0, \ \|e^{t_0\mathcal{H}^{1/2\mu}}f\|_{L^2}<+\infty.
\end{multline}
This characterization proves that there is a regularizing effect for the solutions to the Cauchy problem (\ref{gel1}) in the symmetric Gelfand-Shilov space  $S_{1/2}^{1/2}(\rr^d)$ for any positive time, whereas the smoothing effect for the solutions to the Cauchy problem defined by the evolution equation associated to the fractional harmonic oscillator
\begin{equation}\label{gel2}
\begin{cases}
\partial_tf+\mathcal{H}^sf=0,\\
f|_{t=0}=f_0 \in L^2(\rr^d),
\end{cases}
\end{equation}
with $0<s<1$, occurs for any positive time in the symmetric Gelfand-Shilov space $S_{1/2s}^{1/2s}(\rr^d)$.

In the present work, we investigate the regularizing properties of the Boltzmann equation. More specifically, we study the smoothing effect for the Cauchy problem associated to the radially symmetric spatially homogeneous non-cutoff Boltzmann equation with Maxwellian molecules. In a previous work~\cite{LMPX1}, the linearized radially symmetric non-cutoff Boltzmann operator was shown to behave essentially as a fractional harmonic oscillator $\mathcal{H}^s$, with $0<s<1$. We can directly deduce from the results obtained in~\cite{LMPX1} that there is a regularizing effect for the solutions to the Cauchy problem associated to the linear radially symmetric spatially homogeneous non-cutoff Boltzmann equation with Maxwellian molecules
\begin{equation}\label{gel3}
\begin{cases}
\partial_tf+\mathscr{L}f=0,\\
f|_{t=0}=f_0 \in L^2(\rr^d) \textrm{ radial},
\end{cases}
\end{equation}
in the symmetric Gelfand-Shilov space $S_{1/2s}^{1/2s}(\rr^d)$ for any positive time. This smoothing result is sharp. In this work, we prove that the same smoothing effect in the Gelfand-Shilov space $S_{1/2s}^{1/2s}(\rr^d)$ holds as well for the non-linear equation.

\subsection{The Boltzmann equation}
The Boltzmann equation describes the behaviour of a dilute gas when the only interactions taken into account are binary collisions~\cite{17}. It reads as the equation
\begin{equation}\label{e1}
\begin{cases}
\partial_tf+v\cdot\nabla_{x}f=Q(f,f),\\
f|_{t=0}=f_0,
\end{cases}
\end{equation}
for the density distribution of the particles $f=f(t,x,v) \geq 0$  at time $t$, having position $x \in \rr^d$ and velocity $v \in \rr^d$. The Boltzmann equation derived in 1872 is one of the fundamental equations in mathematical physics and, in particular, a cornerstone of statistical physics.

The term $Q(f,f)$ is the so-called Boltzmann collision operator associated to the Boltzmann bilinear operator
\begin{equation}\label{eq1}
Q(g, f)=\int_{\rr^d}\int_{\SSS^{d-1}}B(v-v_{*},\sigma)(g'_* f'-g_*f)d\sigma dv_*,
\end{equation}
with $d \geq 2$,
where we are using the standard shorthand $f'_*=f(t,x,v'_*)$, $f'=f(t,x,v')$, $f_*=f(t,x,v_*)$, $f=f(t,x,v)$. In this expression, $v, v_*$ and $v',v_*'$ are the velocities in $\rr^d$ of a pair of particles respectively before and after the collision. They are connected through the formulas
$$v'=\frac{v+v_*}{2}+\frac{|v-v_*|}{2}\sigma,\quad   v_*'=\frac{v+v_*}{2}-\frac{|v-v_*|}{2}\sigma,$$
where the parameter $\sigma\in\SSS^{d-1}$ belongs to the unit sphere. Those relations correspond physically to elastic collisions with the conservations of momentum and kinetic energy in the binary collisions
$$\quad v+v_{\ast}=v'+v_{\ast}', \quad |v|^2+|v_{\ast}|^2=|v'|^2+|v_{\ast}'|^2,$$
where $|\cdot|$ is the Euclidean norm on $\rr^d$. The Boltzmann equation is said to be spatially homogeneous when the density distribution of the particles does not depend on the position variable
\begin{equation}\label{vkk1}
\begin{cases}
\partial_tf=Q(f,f),\\
f|_{t=0}=f_0.
\end{cases}
\end{equation}

For monatomic gas, the cross section $B(v-v_*,\sigma)$ is a non-negative function which only depends on the relative velocity $|v-v_*|$ and on the deviation angle $\theta$ defined through the scalar product in $\rr^d$,
$$\cos \theta=k \cdot \sigma, \quad k=\frac{v-v_*}{|v-v_*|}.$$
Without loss of generality, we may assume that $B(v-v_*,\sigma)$ is supported on the set where
$$k \cdot \sigma \geq 0,$$
i.e. where $0 \leq \theta \leq \frac{\pi}{2}$. Otherwise, we can reduce to this situation by the customary symmetrization
$$\tilde{B}(v-v_{*},\sigma)=\big[B(v-v_{*},\sigma)+B(v-v_{*},-\sigma)\big] \un_{\{\sigma \cdot k \geq 0\}},$$
with $\un_A$ being the characteristic function of the set $A$, since the term $f'f_*'$ appearing in the Boltzmann operator $Q(f,f)$ is invariant under the mapping $\sigma \rightarrow -\sigma$.
More specifically, we consider cross sections of  the type
\begin{equation}\label{eq1.01}
B(v-v_*,\sigma)=\Phi(|v-v_*|)b\Big(\frac{v-v_*}{|v-v_*|} \cdot \sigma\Big),
\end{equation}
with a kinetic factor
\begin{equation}\label{sa0}
\Phi(|v-v_*|)=|v-v_*|^{\gamma}, \quad  \gamma \in ]-d,+\infty[.
\end{equation}
The molecules are said to be Maxwellian when the parameter $\gamma=0$.
The second term appearing in the cross sections is a factor related to the deviation angle with a singularity
\begin{equation}\label{sa1}
(\sin \theta)^{d-2}b(\cos \theta)  \substack{\\ \\ \approx \\ \theta \to 0_{+} }  \theta^{-1-2s},
\end{equation}
for\footnote{The notation $a\approx b$ means $a/b$
is bounded from above and below by fixed positive constants.} some  $0 < s <1$. Notice that this singularity is not integrable
$$\int_0^{\frac{\pi}{2}}(\sin \theta)^{d-2}b(\cos \theta)d\theta=+\infty.$$
This non-integrability plays a major role regarding the qualitative behaviour of the solutions of the Boltzmann equation and this feature is essential for the smoothing effect to be present. Indeed, as first observed by Desvillettes for the Kac equation~\cite{D95}, grazing collisions that account for the non-integrability of the angular factor near $\theta=0$
do induce smoothing effects for the solutions of the non-cutoff Kac equation, or more generally for the solutions of the non-cutoff Boltzmann equation. On the other hand, these solutions are at most as regular as the initial data, see e.g. \cite{36}, when the cross section is assumed to be integrable, or after removing the singularity by using a cutoff function (Grad's angular cutoff assumption).

The physical motivation for considering this specific structure of cross sections is derived from particles interacting according to a spherical intermolecular repulsive potential of the form
$$\phi(\rho)=\frac{1}{\rho^{r}}, \quad r>1,$$
with $\rho$ being the distance between two interacting particles. In the physical 3-dimensional space~$\rr^3$, the cross section satisfies the above assumptions with
$$s=\frac{1}{r} \in ]0, 1[, \quad \gamma=1-4s \in ]-3, 1[.$$
For further details on the physics background and the derivation of the Boltzmann equation, we refer the reader to the extensive expositions \cite{17,villani2}.

\subsection{The linearized Boltzmann operator}
We consider the linearization of the Boltzmann equation
$$f=\mu_d+\sqrt{\mu_d}g,$$
around the Maxwellian equilibrium distribution
\begin{equation}\label{maxwe}
\mu_d(v)=(2\pi)^{-\frac{d}{2}}e^{-\frac{|v|^2}{2}}, \quad v \in \rr^d.
\end{equation}
Since $Q(\mu_d,\mu_d)=0$ by the conservation of the kinetic energy, the Boltzmann operator $Q(f,f)$ can be split into three terms
$$Q(\mu_d+\sqrt{\mu_d}g,\mu_d+\sqrt{\mu_d}g)=Q(\mu_d,\sqrt{\mu_d}g)+Q(\sqrt{\mu_d}g,\mu_d)+Q(\sqrt{\mu_d}g,\sqrt{\mu_d}g),$$
whose linearized part is
$Q(\mu_d,\sqrt{\mu_d}g)+Q(\sqrt{\mu_d}g,\mu_d).$
Setting
\begin{equation}\label{jan1}
\mathscr{L}g=\mathscr{L}_1g+\mathscr{L}_2g,
\end{equation}
with
\begin{equation}\label{jan1.1}
\mathscr{L}_1g=-\mu_d^{-1/2}Q(\mu_d,\mu_d^{1/2}g), \quad \mathscr{L}_2g=-\mu_d^{-1/2}Q(\mu_d^{1/2}g,\mu_d),
\end{equation}
the original spatially homogeneous Boltzmann equation (\ref{vkk1}) is reduced to the Cauchy problem for the fluctuation
\begin{equation}\label{boltz}
\begin{cases}
\partial_tg+\mathscr{L}g=\mu_d^{-1/2}Q(\sqrt{\mu_d}g,\sqrt{\mu_d}g),\\
g|_{t=0}=g_0.
\end{cases}
\end{equation}
The Boltzmann operator is local in the time and position variables and from now on, we consider it as acting only in the velocity variable.
This linearized operator is known \cite{17} to be an unbounded symmetric operator on $L^2(\rr^d_{v})$ (acting in the velocity variable) such that its Dirichlet form satisfies
$$(\mathscr{L}g,g)_{L^2(\rr^d_{v})} \geq 0.$$
Setting
$$\mathbf{P}g=(a+b \cdot v+c|v|^2)\mu_d^{1/2},$$
with $a,c \in \rr$, $b \in \rr^d$, the $L^2$-orthogonal projection onto the space of collisional invariants
\begin{equation}\label{coli}
\mathcal{N}=\textrm{Span}\big\{\mu_d^{1/2},v_1 \mu_d^{1/2},...,v_d\mu_d^{1/2},|v|^2\mu_d^{1/2}\big\},
\end{equation}
we have
\begin{equation}\label{ker}
(\mathscr{L}g,g)_{L^2(\rr^d)}=0 \Leftrightarrow g=\mathbf{P}g.
\end{equation}
It was noticed forty years ago by Cercignani \cite{cercignani}  that the linearized Boltzmann operator with Maxwellian molecules behaves like a fractional diffusive operator.
Over the time, this point of view transformed into the following widespread heuristic conjecture on the diffusive behavior of the Boltzmann operator as a flat fractional Laplacian \cite{al2009,al-1,amuxy-2,pao1,pao2,villani2}:
$$f \mapsto Q(\mu_d,f) \sim -(-\Delta_v)^sf+ \textrm{ lower order terms},$$
with $0<s<1$ being the parameter appearing in the singularity assumption (\ref{sa1}). See \cite{lmp,MoXu,M-X2} for works related to this simplified model of the non-cutoff Boltzmann equation. In the physical 3-dimensional case, we recently unveiled the exact nature of the linearized non-cutoff Boltzmann operator with Maxwellian molecules by proving that this operator is equal to the fractional linearized Landau operator with Maxwellian molecules \cite{LMPX2} (Theorem~2.3),
$$\mathscr{L}=a(\mathcal{H},\Delta_{\SSS^2}) \mathscr{L}_L^s,$$
up to a positive bounded isomorphism on $L^2(\rr^3)$,
$$\exists c>0, \forall f \in L^2(\rr^3), \quad c\|f\|_{L^2}^2 \leq (a(\mathcal{H},\Delta_{\SSS^2})f,f)_{L^2} \leq \frac{1}{c}\|f\|_{L^2}^2,$$
commuting with the harmonic oscillator
$\mathcal H=-\Delta_v+\frac{|v|^2}{4}$
and the Laplace-Beltrami operator
$$\Delta_{\SSS^{2}}=\frac{1}{2}\sum_{\substack{1 \leq j,k \leq 3 \\ j \neq k}}(v_j \partial_k-v_k \partial_j)^2,$$
on the unit sphere $\SSS^{2}$. The operator
$$\mathscr{L}_Lg=-\mu_3^{-1/2}Q_L(\mu_3,\mu_3^{1/2}g)-\mu_3^{-1/2}Q_L(\mu_3^{1/2}g,\mu_3),$$
refers to the linearization of the Landau collision operator around the Maxwellian distribution (\ref{maxwe}) with $d=3$,
$$Q_L(g, f)=\nabla_v \cdot \Big(\int_{\RR^3}a(v-v_*)\big(g(t,x,v_*)(\nabla_v f)(t,x,v)-(\nabla_v g)(t,x,v_*)f(t,x,v)\big)d v_*\Big),$$
where $a=(a_{i,j})_{1 \leq i,j \leq 3}$ stands for the non-negative symmetric matrix
\begin{equation}\label{landau_collision1}
a(v)=(|v|^2\Id -v\otimes v) \in M_3(\rr),
\end{equation}
which is equal to the differential operator \cite{LMPX2} (Proposition~2.1),
$$\mathscr{L}_L=2\mathcal{H}-\Delta_{\SSS^{2}}-3+(-2\mathcal{H}+\Delta_{\SSS^{2}}+3)\mathbb{P}_1 +(-2\mathcal{H}-\Delta_{\SSS^{2}}+3)\mathbb{P}_2,$$
when acting on the Schwartz space $\mathscr{S}(\rr^3)$, where $\mathbb{P}_{k}$ are the orthogonal projections onto the Hermite basis defined in Section~\ref{6.sec.harmo}. More generally, the linearized non-cutoff Boltzmann operator with general molecules $\gamma \in ]-3,+\infty[$ was shown to satisfy the following coercive estimates \cite{LMPX2} (Theorem~2.5),
\begin{multline*}
\forall f \in \mathscr{S}(\rr^3), \quad   \|\mathcal{H}^{\frac{s}{2}}\langle v \rangle^{\frac{\gamma}{2}}(1-{\bf P})f\|_{L^2}^2 +\|(-\Delta_{\SSS^2})^{\frac{s}{2}}\langle v \rangle^{\frac{\gamma}{2}}(1-{\bf P})f\|_{L^2}^2 \\ \lesssim (\mathscr{L}f,f)_{L^2} \lesssim  \|\mathcal{H}^{\frac{s}{2}}\langle v \rangle^{\frac{\gamma}{2}}(1-{\bf P})f\|_{L^2}^2 +\|(-\Delta_{\SSS^2})^{\frac{s}{2}}\langle v \rangle^{\frac{\gamma}{2}}(1-{\bf P})f\|_{L^2}^2.
\end{multline*}
For Maxwellian molecules, the spectrum of the linearized Boltzmann operator is only composed of eigenvalues explicitly computed in \cite{wang}. See also \cite{bobylev,17,dolera}.
In the physical 3-dimensional case, we define
\begin{equation}\label{new001}
\beta(\theta)=|\SSS^{1}||\sin 2\theta|b(\cos 2\theta).
\end{equation}
This function satisfies 
\begin{equation}\label{new001b}
\beta(\theta) \substack{\\ \\ \approx \\ \theta \to 0 } \val\theta^{-1-2s},
\end{equation}
when the singularity assumption (\ref{sa1}) holds.
For $\sigma=(\cos \phi \sin \alpha,\sin \phi \sin \alpha,\cos \alpha) \in \SSS^2,$
with $\alpha \in [0,\pi]$ and $\phi \in [0,2\pi)$, the real spherical harmonics $Y_l^m(\sigma)$ with $l \geq 0$, $-l \leq m \leq l$, are defined as $Y_0^0(\sigma)=(4\pi)^{-1/2}$ and for any $l \geq 1$,
$$Y_l^m(\sigma)=\begin{cases} \displaystyle \sqrt{\frac{2l+1}{4\pi}}P_l(\cos \alpha), & \mbox{if } m=0 \\
\displaystyle \sqrt{\frac{2l+1}{2\pi}\frac{(l-m)!}{(l+m)!}}P_l^m(\cos \alpha)\cos m \phi  & \mbox{if } m=1,...,l \\
\displaystyle \sqrt{\frac{2l+1}{2\pi}\frac{(l+m)!}{(l-m)!}}P_l^{-m}(\cos \alpha)\sin m \phi  & \mbox{if } m=-l,...,-1,
\end{cases}$$
where $P_l$ stands for the $l$-th Legendre polynomial and $P_l^m$ the associated Legendre functions of the first kind of order $l$ and degree $m$. The family $(Y_l^m)_{l \geq 0, -l \leq m \leq l}$ constitutes an orthonormal basis of the space $L^2(\SSS^2,d\sigma)$ with $d\sigma$ being the surface measure on $\SSS^2$. We set for any $n,l \geq 0$, $-l \leq m \leq l$,
\begin{equation}\label{eig11}
\varphi_{n,l,m}(v)=2^{-1/4}\sqrt{\frac{n!}{\Gamma(n+l+\frac{3}{2})}}\Big(\frac{|v|}{\sqrt{2}}\Big)^lL_n^{[l+\frac{1}{2}]}\Big(\frac{|v|^2}{2}\Big)e^{-\frac{|v|^2}{4}}Y_l^m\Big(\frac{v}{|v|}\Big),
\end{equation}
where $L_n^{[l+\frac{1}{2}]}$ are the generalized Laguerre polynomials. The family $(\varphi_{n,l,m})_{n,l \geq 0, |m| \leq l}$ is an orthonormal basis of $L^2(\rr^3)$
composed of eigenvectors of the harmonic oscillator and the Laplace-Beltrami operator on the unit sphere $\SSS^2$,
\begin{equation}\label{k4.1}
\Big(-\Delta_v+\frac{|v|^2}{4}-\frac{3}{2}\Big)\varphi_{n,l,m}=(2n+l)\varphi_{n,l,m}, \quad -\Delta_{\SSS^2}\varphi_{n,l,m}=l(l+1)\varphi_{n,l,m}.
\end{equation}
The linearized non-cutoff Boltzmann operator with Maxwellian molecules is also diagonal in this orthonormal basis $(\varphi_{n,l,m})_{n,l \geq 0, |m| \leq l}$. In the cutoff case, i.e., when
$b(\cos \theta)\sin \theta \in L^1([0,\pi/2]),$
it was shown in \cite{wang} that
\begin{equation}\label{k4.23}
\mathscr{L}\varphi_{n,l,m}=\lambda_{B}(n,l,m)\varphi_{n,l,m}, \quad n,l \geq 0, \ -l \leq m \leq l,
\end{equation}
with
\begin{multline}\label{k4.24}
\lambda_B(n,l,m)= \int_{-\frac{\pi}{4}}^{\frac{\pi}{4}}\beta(\theta)\big(1+\delta_{n,0}\delta_{l,0}-P_l(\cos \theta)(\cos \theta)^{2n+l}-P_l(\sin \theta)(\sin \theta)^{2n+l}\big)d\theta,
\end{multline}
where $P_l$ are the Legendre polynomials
 \begin{equation}\label{rodrigues}
 P_l(x)=\frac{1}{2^l l!}\frac{d^l}{dx^l}(x^2-1)^l, \quad l \geq 0.
\end{equation}
This diagonalization of the linearized Boltzmann operator with Maxwellian molecules holds as well in the non-cutoff case, see e.g. \cite{bobylev,17,dolera,LMPX2}. The space of the collisional invariants (\ref{coli}) may be expressed throughout this diagonalizing basis as
$$\mathcal{N}= \textrm{Span}\big\{\varphi_{0,0,0},\varphi_{0,1,-1},\varphi_{0,1,0},\varphi_{0,1,1},\varphi_{1,0,0}\big\}.$$
Notice that all the scalar products
\begin{equation}\label{rumi10}
(f,\varphi_{n,l,m})_{L^2}=0, \quad n,l \geq 0,\ -l \leq m \leq l,\ (l,m) \neq (0,0),
\end{equation}
are zero when $f \in L^2(\rr^3)$ is a radial function.
When acting on radial functions, the linearized non-cutoff Boltzmann operator with Maxwellian molecules is therefore equal to
\begin{equation}\label{rr2}
\mathscr{L}=\sum_{k=1}^{+\infty}\lambda_{2k}\mathbb{P}_{2k},
\end{equation}
with
\begin{equation}\label{tea1}
\lambda_2=0, \qquad
\lambda_{2k}=\int^{\frac{\pi}{4}}_{-\frac{\pi}{4}}\beta(\theta)(1-(\cos\theta)^{2k}-(\sin \theta)^{2k})d\theta \geq 0, \quad  k \geq 2,
\end{equation}
where $\mathbb{P}_{2k}$ stands for the orthogonal Hermite projections onto the energy level $2k$ described in Section \ref{6.sec.harmo}, since
$$\mathbb{P}_{2k}f=(f,\varphi_{k,0,0})_{L^2}\varphi_{k,0,0},$$
when $f \in L^2(\rr^3)$ is a radial function.
The assumption (\ref{new001b}) on the cross section $\beta$ implies the following asymptotics
\begin{equation}\label{ekk1}
\lambda_{2k} \sim \int^{\frac{\pi}{4}}_{-\frac{\pi}{4}}\beta(\theta)(1-(\cos\theta)^{2k})d\theta \approx k^s,
\end{equation}
when $k \rightarrow +\infty$, see \cite{LMPX1} (Corollary~2.7).

\section{Statement of the main result}\label{statements}

We consider the linearized non-cutoff Boltzmann operator with Maxwellian molecules
$$\mathscr{L}f=-\mu_3^{-1/2}Q(\mu_3,\mu_3^{1/2}f)-\mu_3^{-1/2}Q(\mu_3^{1/2}f,\mu_3),$$
acting on the radially symmetric Schwartz space on $\rr^3$ (see Section \ref{6.sec.radia}),
{\small
\begin{equation}\label{2.radia}
\mathscr{S}_{r}(\rr^3)=\big\{f\in\mathscr{S}(\rr^3),\ \forall v\in\rr^3,\forall A\in O(3),\ f(v)=f(Av)\big\}=\big\{f(\val v)\big\}_{\substack{f \text{ even }\\ f \in \mathscr S(\R)}},
\end{equation}
}\noindent
where $O(3)$ stands for the orthogonal group of $\rr^3$. We recall that the case of Maxwellian molecules corresponds to the case when $\gamma=0$ in the kinetic factor (\ref{sa0}) and that the non-negative cross section $b(\cos \theta)$ is assumed to be supported where $\cos \theta \geq 0$ and to satisfy the assumption (\ref{sa1}).

Following the Bobylev's theory~\cite{bobylev}, we solve explicitly the Cauchy problem associated to the non-cutoff radially symmetric spatially homogeneous Boltzmann equation with Maxwellian molecules for small initial radial $L^2$-fluctuations around the standard Maxwellian distribution.
In~\cite{bobylev} (p. 215), Bobylev constructs explicit global radial solutions for initial radial $L^2$-fluctuations
$$f_0=\mu_3+\sqrt{\mu_3}g_0, \qquad g_0=\sum_{n=2}^{+\infty}b_{n}(0)\varphi_{n,0,0},$$
satisfying
\begin{equation}\label{rrr2}
\sup_{n \geq 2}\sqrt[n]{|b_n(0)|\sqrt{\frac{\pi^{1/2}\ n!}{2\Gamma(n+\frac{3}{2})}}}<\frac{3}{7},
\end{equation}
and establishes the exponential return to equilibrium for the density distribution of the particles
$$f=\mu+\sqrt{\mu}g,$$
in the $L^{\infty}(\rr_v^3)$-norm
$$\exists C>0, \forall t \geq 0, \quad \|f(t)-\mu\|_{L^{\infty}} \leq Ce^{-\lambda_4t}.$$
In this work, we do not request the specific structure (\ref{rrr2}) and perform the construction of explicit global radial solutions for any sufficiently small initial radial $L^2$-fluctuation. The main novelty of the present work relates to the property of exponential convergence to zero for the fluctuation which is established in a specific weighted space emphasizing that the non-cutoff radially symmetric spatially homogeneous Boltzmann equation enjoys regularizing properties in the Gelfand-Shilov space $S_{1/2s}^{1/2s}(\rr^3)$ for any positive time.

Regarding the smoothing features of the Boltzmann equation, we recall that the non-cutoff spatially homogeneous Boltzmann equation is known to enjoy
a $\mathscr{S}(\rr^d)$-regularizing effect for the weak solutions to the Cauchy problem~\cite{DW2004}. Regarding the Gevrey regularity, Ukai showed in~\cite{ukai} that the Cauchy problem for the Boltzmann equation has a unique local solution in Gevrey classes.  Then, Desvillettes, Furioli and Terraneo proved in~\cite{desv-fl-terr} the propagation of Gevrey regularity for solutions of the Boltzmann equation with Maxwellian molecules. For mild singularities $\gamma+2s<1$, Morimoto and Ukai proved in~\cite{ukai1} the $G^{1/2s}$-Gevrey regularity of smooth Maxwellian decay solutions to the Cauchy problem of the spatially homogeneous Boltzmann equation with a modified kinetic factor $\Phi(|v-v_*|)=\langle v-v_* \rangle^{\gamma}$. This result for mild singularities was recently extended by Zhang and Yin~\cite{zhang} for the standard kinetic factor $\Phi(|v-v_*|)=|v-v_*|^{\gamma}$. In~\cite{MUXY-DCDS}, Morimoto, Ukai, Xu and Yang have established the property of $G^{1/s}$-Gevrey smoothing effect for the weak solutions to the Cauchy problem associated to the linearized spatially homogeneous Boltzmann equation with Maxwellian molecules when $0<s<1$. On the other hand, Lekrine and Xu proved in~\cite{LX09} the property of $G^{1/2s'}$-Gevrey smoothing effect for the weak solutions to the Cauchy problem associated to the radially symmetric spatially homogeneous Boltzmann equation with Maxwellian molecules for any $0<s'<s$, when the singularity is mild $0<s<1/2$. This result was then completed by Glangetas and Najeme who established in~\cite{GN} the analytic smoothing effect in the case when $1/2<s<1$.

Setting
\begin{equation}\label{ar1}
\alpha_{2n, 2m}=\sqrt{C_{2n+2m}^{2n}}\Big(\int_{-\frac\pi{4}}^{\frac\pi{4}}\beta(\theta)(\sin \theta)^{2n} (\cos \theta)^{2m}d\theta\Big), \quad n\geq 1, \ m \geq 0,
\end{equation}
\begin{equation}\label{ar2}
\alpha_{0, 2m}=-\Big(\int_{-\frac\pi{4}}^{\frac\pi{4}}\beta(\theta)(1-(\cos \theta)^{2m})d\theta\Big), \quad m \geq 1, \quad \alpha_{0,0}=0,
\end{equation}
where $C_n^k=\frac{n!}{k!(n-k)!}$ stands for the binomials coefficients, we consider the infinite system of differential equations
\begin{equation}\label{evv3}
\left\{
\begin{array}{c}
    \partial_t b_0(t) = 0 \\
 \displaystyle   \forall n \geq 1, \quad \partial_t b_{n}(t)+\lambda_{2n} b_{n}(t) =  \displaystyle
\alpha_{0,2n}b_0(t)b_{n}(t)\\
\displaystyle \hspace{2cm} +\sum_{\substack{k+l=n \\ k \geq 1, l \geq 0}}\alpha_{2k,2l}\sqrt{\frac{(2k+2l+1)}{(2k+1)(2l+1)}}b_{k}(t)b_{l}(t),
\end{array}
\right.
\end{equation}
where $\lambda_{2n}$ stands for the eigenvalue of the linearized radially symmetric Boltzmann operator (\ref{tea1}).
This system is  triangular
\begin{equation}\label{evv4}
\left\{
\begin{array}{c}
\displaystyle \forall t \geq 0, \quad b_0(t)=b_0(0),\\
\displaystyle \forall t \geq 0, \quad b_1(t)=b_1(0),\\
\displaystyle \forall n \geq 2, \forall t \geq 0, \quad \partial_t b_{n}(t)+\lambda_{2n}(1+b_0(0)) b_{n}(t)\\
\displaystyle \hspace{2cm} =\sum_{\substack{k+l=n \\ k \geq 1, l \geq 1}}\alpha_{2k,2l}\sqrt{\frac{(2k+2l+1)}{(2k+1)(2l+1)}}b_{k}(t)b_{l}(t),
\end{array}
\right.
\end{equation}
since the $(n+1)^{\textrm{th}}$ equation is a linear differential equation for the function $b_{n}$ with a right-hand-side involving only the functions $(b_{k})_{1 \leq k \leq n-1}$. This system may therefore be explicitly solved while solving a sequence of linear differential equations. This allows to solve explicitly the non-cutoff radially symmetric spatially homogeneous Boltzmann equation:

\bigskip

\begin{theorem}\label{thnonlinear11}
Let $0 <\delta < 1$ be a positive constant.
There exists a positive constant $\eps_0>0$ such that if $g_0 \in \mathcal{N}^{\perp}$ is a radial  $L^2(\rr^3)$-function satisfying $\|g_0\|_{L^2} \leq \eps_0$, then the Cauchy problem for the fluctuation associated to the non-cutoff spatially homogeneous Boltzmann equation with Maxwellian molecules
$$\begin{cases}
\partial_tg+\mathscr{L}g=\mu_3^{-1/2}Q(\sqrt{\mu_3}g,\sqrt{\mu_3}g),\\
g|_{t=0}=g_0,
\end{cases}$$
has a unique global radial solution $g \in L^{\infty}(\rr_t^+,L^2(\rr_v^3))$ given by
$$g(t)=\sum_{n=0}^{+\infty}b_{n}(t)\varphi_{n,0,0},$$
where $(\varphi_{n,0,0})_{n\in \N}$ are the functions defined in \emph{(\ref{eig11})} and where the functions $(b_{n}(t))_{n \geq 0}$ are the solutions of the system of differential equations $(\ref{evv4})$ with initial conditions
$$\forall n \geq 0, \ b_{n}(t)|_{t=0}=(g_0,\varphi_{n,0,0})_{L^2}.$$
Furthermore, this fluctuation around the Maxwellian distribution is exponentially convergent to zero in the following weighted $L^2$-space
\begin{equation}\label{rr1}
\forall t \geq 0, \quad \|e^{\frac{t}{2}\mathscr{L}}g(t)\|_{L^2}=\Big(\sum_{n=0}^{+\infty}e^{\lambda_{2n}t}|b_{n}(t)|^2\Big)^{1/2} \leq e^{-\frac{\lambda_4}{2}(1-\delta) t}\|g_0\|_{L^2}
\end{equation}
and belongs to the Gelfand-Shilov class $S_{1/2s}^{1/2s}(\rr^3)$ for any positive time
\begin{equation}\label{rrr1}
\forall t >0, \quad g(t) \in S_{1/2s}^{1/2s}(\rr^3),
\end{equation}
where $0<s<1$ is the parameter appearing in the singularity assumption \emph{(\ref{sa1})}.

\end{theorem}

\bigskip

This result emphasizes that the non-cutoff radially symmetric spatially homogeneous Boltzmann equation enjoys specific Gelfand-Shilov regularizing properties depending directly on the value of the parameter $0<s<1$ in the singularity assumption (\ref{sa1}).  
The Gelfand-Shilov smoothing effect
$$\forall t >0, \quad g(t) \in S_{1/2s}^{1/2s}(\rr^3),$$
shows that the fluctuation enjoys a $G^{1/2s}$-Gevrey smoothing effect for any positive time $t >0$,
$$g(t) \in G^{1/2s}(\rr^3),$$
together with an additional decay at infinity exactly controlled by an exponential term of the type $e^{-c_t|v|^{2s}}$, for some $c_t>0$. In particular, this result points out an ultra-analytic smoothing effect for the range of parameter $1/2<s<1$, since in this case the solution to the non-cutoff radially symmetric spatially homogeneous Boltzmann equation is smoother then analytic
$$G^{1/2s} \subset G^1,$$
when $1/2<s<1$. Up to our knowledge, this is the first result of ultra-analytic smoothing obtained for the Boltzmann equation since in previous results in the literature recalled at the beginning of this section, the $G^{1/2s}$-Gevrey regularity was only obtained in~\cite{ukai1,zhang} for smooth Maxwellian decay solutions in the case of mild singularities $0<s<1/2$. In addition, we obtain in Theorem~\ref{thnonlinear11} a precise control of the smoothing effect (\ref{rr1}) implying the return to equilibrium in a weighted $L^2$-space.

The Gelfand-Shilov smoothing effect
$$\forall t >0, \quad g(t) \in S_{1/2s}^{1/2s}(\rr^3),$$
is a direct consequence of the a priori estimate (\ref{rr1}) and the spectral asymptotics (\ref{ekk1}). This result implies in particular that there exists a positive constant $c>0$ such that any global radial solution $g \in L^{\infty}(\rr_t^+,L^2(\rr_v^3))$ associated to an initial radial fluctuation $g_0 \in \mathcal{N}^{\perp}$ satisfying $\|g_0\|_{L^2} \leq \eps_0$, given by Theorem~\ref{thnonlinear11} satisfies  
the estimate
\begin{equation}\label{zok}
\forall t \geq 0, \quad \|e^{ct\mathcal{H}^r}g(t)\|_{L^2} \leq \|g_0\|_{L^2}.
\end{equation}
with $r=s$ the parameter appearing in the singularity assumption (\ref{sa1}). The estimate (\ref{zok}) is sharp with respect to the index $r$. We can check this sharpness by solving explicitly the triangular system of differential equations (\ref{evv4}). Indeed, when the initial fluctuation is a radial function satisfying $g_0 \in \mathcal{N}^{\perp}$, i.e. $b_0(0)=b_1(0)=0$, we deduce from Lemma~\ref{lemk10} that non-linear effects do not appear before the component $b_4$,
\begin{equation}\label{kek5}
\forall t \geq 0, \quad b_0(t)=b_1(t)=0,
\end{equation}
\begin{equation}\label{kek6}
\forall t \geq 0, \quad b_2(t)=b_2(0)e^{-\lambda_4t}, \quad b_3(t)=b_3(0)e^{-\lambda_6t},
\end{equation}
\begin{multline*}
\forall t \geq 0, \quad b_4(t)=\Big[b_4(0)-\frac{3}{5}\sqrt{C_8^4}\frac{b_2(0)^2}{\lambda_8-2\lambda_4}\Big(\int_{-\frac{\pi}{4}}^{\frac{\pi}{4}}\beta(\theta)(\sin \theta)^4(\cos \theta)^4 d\theta\Big)\Big]e^{-\lambda_8 t}\\
+\frac{3}{5}\sqrt{C_8^4}\frac{b_2(0)^2}{\lambda_8-2\lambda_4}\Big(\int_{-\frac{\pi}{4}}^{\frac{\pi}{4}}\beta(\theta)(\sin \theta)^4(\cos \theta)^4 d\theta\Big)e^{-2\lambda_4t},
\end{multline*}
and we establish by induction that for any $n \geq 1$, there exist some constants $\gamma_{j_1,j_2,...,j_k}$ such that
\begin{equation}\label{ek1111}
\forall t \geq 0, \quad b_{n}(t)=\sum_{1 \leq k \leq n}\sum_{\substack{j_1+j_2+...+j_k=n\\ j_1,...,j_k \geq 2}}\gamma_{j_1,j_2,...,j_k}e^{-(\lambda_{2j_1}+\lambda_{2j_2}+...+\lambda_{2j_k})t}.
\end{equation}
According to Lemma~\ref{lemk10}, we notice that
$$\lambda_{2n} \leq \lambda_{2j_1}+\lambda_{2j_2}+...+\lambda_{2j_k}.$$
Furthermore, the choice of particular initial radial fluctuations satisfying 
$$b_n(0)=(g_0,\varphi_{n,0,0})_{L^2}=\eps_0, \quad \forall k \neq n, \quad b_k(0)=(g_0,\varphi_{k,0,0})_{L^2}=0,$$
for some $n \geq 2$, implies that 
$$\forall t\geq 0, \quad b_n(t)=b_n(0)e^{-\lambda_{2n}t}.$$
We therefore deduce from (\ref{ekk1}) that the uniform control (\ref{zok}) cannot hold if $r>s$, since the term
$$\|e^{ct\mathcal{H}^r}g_n(t)\|_{L^2} \geq \eps_0 e^{t(c(2n+\frac{3}{2})^r-\lambda_{2n})},$$
tends to infinity for any fixed $t>0$, when $n \to +\infty$.

\section{Proof of the main result}

\subsection{Preliminary results}
This first section is dedicated to recall some elements of the Bobylev's theory~\cite{bobylev}, which will be needed for our analysis. For the sake of completeness and the reader's convenience, we provide complete proofs of all these results.

We begin this section by emphasizing some properties of no resonances satisfied by the eigenvalues of the linearized radially symmetric Boltzmann operator (\ref{tea1}), which are playing a basic role in the present analysis:

\bigskip

\begin{lemma}\label{lemk10}
The eigenvalues of the linearized radially symmetric Boltzmann operator
$$\lambda_{2n}=\int^{\frac{\pi}{4}}_{-\frac{\pi}{4}}\beta(\theta)(1-(\cos\theta)^{2n}-(\sin \theta)^{2n})d\theta \geq 0, \quad n \geq 2,$$
satisfy the following estimates
$$\forall j,k \geq 2, \quad \lambda_{2j+2k} <\lambda_{2j}+\lambda_{2k}.$$
\end{lemma}

\bigskip

\begin{proof}
We first prove the estimate $\lambda_8 < 2 \lambda_4.$
This is a direct consequence of the inequality
$$1-\cos ^8\theta-\sin^8\theta<2(1-\cos^4\theta-\sin^4\theta), \quad \theta \in I=\Big[-\frac{\pi}{4},\frac{\pi}{4}\Big] \setminus \{0\},$$
which holds true since
$$2(1-\cos^4\theta-\sin^4\theta)-1+\cos ^8\theta+\sin^8\theta=2\cos^4\theta (\cos^2\theta-1)^2=2\cos^4\theta \sin^4\theta>0.$$
Then, we prove that
$$\lambda_{2j+2k}<\lambda_{2j}+\lambda_{2k},\quad j,k\ge2, \ (j,k) \neq (2,2).$$
This is a direct consequence of the inequality
$$1+\cos^{2j+2k}\theta+\sin^{2j+2k}\theta -\cos^{2j}\theta-\sin^{2j}\theta-\cos^{2k}\theta-\sin^{2k}\theta>0, \quad \theta \in I.$$
To that end, we check that for any $j,k\ge 2$, $(j,k)\not=(2,2)$,
\begin{equation}\label{trivineq0}
1+(1-t)^{j+k}+t^{j+k}-(1-t)^j-t^j-(1-t)^k-t^k>0, \quad  t\in (0,1/2].
\end{equation}
The latter is equivalent to
\begin{equation}\label{trivineq}
(1-(1-t)^k)(1-(1-t)^j)> t^j+t^k-t^{j+k}, \quad  t\in (0,1/2].
\end{equation}
Since
\begin{equation}\label{vkk2}
1-(1-t)^j=t\sum_{l=0}^{j-1}(1-t)^{l},
\end{equation}
it follows that for any $j,k\ge 2$, $(j,k)\not=(2,2)$,
$$(1-(1-t)^k)(1-(1-t)^j)=t^2\sum_{0\le l_1 \leq j-1\atop 0\le l_2 \leq k-1}(1-t)^{l_1+l_2}\ge t^2(1-t)^{j-2}+t^2(1-t)^{k-2}.$$
When $0 < t\le 1/2$, the estimate $1-t\ge t$ implies that
$$(1-(1-t)^k)(1-(1-t)^j)\ge t^j+t^k > t^j+t^k-t^{j+k}.$$
This proves that the estimate \eqref{trivineq} holds for $0 < t\le 1/2$.
\end{proof}

\bigskip

\noindent
The algebraic identities in the following lemma are instrumental in the proof of the Lemma~\ref{proposition111}:

\bigskip

\begin{lemma}\label{eigbolt}
Let $\varphi_{n,l,m} \in L^2(\rr^3)$ be the functions defined in \emph{(\ref{eig11})} and $\psi_n \in L^2(\rr)$ the Hermite functions defined in Section~\ref{6.sec.harmo}. Then, we have for all $n \geq 0$,
$$(-1)^n \sqrt{2n+1}\widehat{\mu_3^{1/2}\varphi_{n,0,0}}(\xi)=\widehat{\mu_1^{1/2}\psi_{2n}}(|\xi|), \quad \xi \in \rr^3.$$
\end{lemma}

\bigskip

\begin{proof}
We may write for any $n \geq 0$,
\begin{multline*}
\mu_1^{1/2}(v)\psi_n(v)=\frac{1}{(4\pi)^{\frac{1}{4}}}\phi_n\Big(\frac{v}{\sqrt{2}}\Big)e^{-\frac{v^2}{4}}=\frac{(-1)^n}{\sqrt{2^{n+1}\pi n! }} \Big[\frac{d^n}{dx^n}(e^{-x^2})\Big]\Big|_{x=\frac{v}{\sqrt{2}}}\\
=\frac{(-1)^ni^n}{\sqrt{2\pi n! }}D_v^n\big(e^{-\frac{v^2}{2}}\big).
\end{multline*}
It follows that
\begin{equation}\label{rumi1}
\widehat{\mu_1^{1/2}\psi_n}(\xi)=\frac{(-1)^ni^n}{\sqrt{2\pi n! }}\xi^n\widehat{e^{-\frac{v^2}{2}}}(\xi)=\frac{(-1)^ni^n}{\sqrt{n!}}\xi^ne^{-\frac{\xi^2}{2}},
\end{equation}
since
\begin{equation}\label{eq4}
\widehat{\big(e^{-\frac{\alpha}{2} \val v^2}\big)}(\xi)=\int_{\R^{d}}e^{-\frac{\alpha}{2} \val v^2}e^{-iv\cdot\xi}dv=\frac{(2\pi)^{\frac{d}{2}}}{\alpha^{\frac{d}{2}}}e^{-\frac{\val\xi^2}{2\alpha}},
\end{equation}
when $\alpha>0$.
According to (\ref{rumi1}), it is sufficient to check that
\begin{equation}\label{rumi2}
\widehat{\mu_3^{1/2}\varphi_{n,0,0}}(\xi)=\frac{1}{\sqrt{(2n+1)!}}|\xi|^{2n}e^{-\frac{|\xi|^2}{2}}.
\end{equation}
By using the standard formula, see e.g. (6.2.15) in~\cite{askey} (Chapter 6),
$$L_n^{[\frac{1}{2}]}(x)=\frac{e^xx^{-\frac{1}{4}}}{n!}\int_0^{+\infty}t^{n+\frac{1}{4}}J_{1/2}(2\sqrt{xt})e^{-t}dt, \quad x>0,$$
where $J_{1/2}$ stands for the Bessel function
\begin{equation}\label{rumi5}
J_{1/2}(x)=\sqrt{\frac{2}{\pi x}}\sin x=\sqrt{\frac{x}{2\pi}}\int_0^{\pi}e^{ix\cos \theta} \sin \theta d\theta, \quad x>0,
\end{equation}
we deduce from (\ref{eig11}) that
\begin{align*}
\varphi_{n,0,0}(v)=& \ 2^{-5/4}\sqrt{\frac{n!}{\pi\Gamma(n+\frac{3}{2})}}L_n^{[\frac{1}{2}]}\Big(\frac{|v|^2}{2}\Big)e^{-\frac{|v|^2}{4}}\\
= & \ \frac{1}{2\sqrt{\pi|v|\Gamma(n+\frac{3}{2})n!}}e^{\frac{|v|^2}{4}}\int_0^{+\infty}t^{n+\frac{1}{4}}J_{1/2}(|v|\sqrt{2t})e^{-t}dt.
\end{align*}
It follows from (\ref{maxwe}) that
\begin{equation}\label{rumi3}
\mu_3^{1/2}\varphi_{n,0,0}(v)= \frac{1}{2^{n+2}\pi^{\frac{5}{4}}\sqrt{|v|\Gamma(n+\frac{3}{2})n!}}\int_0^{+\infty}r^{2n+\frac{3}{2}}J_{1/2}(r|v|)e^{-\frac{r^2}{2}}dr.
\end{equation}
On the other hand, we may compute the inverse Fourier transform
\begin{multline*}
F(v)=\frac{1}{(2\pi)^3}\int_{\rr^3}\frac{1}{\sqrt{(2n+1)!}}|\xi|^{2n}e^{-\frac{|\xi|^2}{2}}e^{iv \cdot \xi}d\xi\\ =\frac{1}{(2\pi)^3\sqrt{(2n+1)!}}\int_{r=0}^{+\infty}\int_{\theta=0}^{\pi}\int_{\phi=0}^{2\pi}r^{2n+2}e^{-\frac{r^2}{2}}e^{i|v|r \cos \theta}\sin \theta drd\theta d\phi.
\end{multline*}
It follows from (\ref{rumi5}) that
\begin{equation}\label{rumi6}
F(v)=\frac{1}{(2\pi)^{\frac{3}{2}}\sqrt{|v|(2n+1)!}}\int_{0}^{+\infty}r^{2n+\frac{3}{2}}J_{1/2}(r|v|)e^{-\frac{r^2}{2}}dr.
\end{equation}
By using that
$$\Gamma\Big(n+\frac{3}{2}\Big)=\Big(n+\frac{1}{2}\Big)\times \Big(n-\frac{1}{2}\Big)\times ...\times \frac{3}{2}\times \frac{1}{2}\times \Gamma\Big(\frac{1}{2}\Big)=\frac{(2n+1)!!\sqrt{\pi}}{2^{n+1}}=\frac{(2n+1)!\sqrt{\pi}}{2^{2n+1}n!},$$
we notice that
$$2^{n+2}\pi^{\frac{5}{4}}\sqrt{\Gamma\Big(n+\frac{3}{2}\Big)n!}=(2\pi)^{\frac{3}{2}}\sqrt{(2n+1)!}.$$
Then, the identity (\ref{rumi2}) follows from (\ref{rumi3}) and (\ref{rumi6}).
\end{proof}

\bigskip

\noindent
The following lemma emphasizes noticeable algebraic identities satisfied by the Boltzmann collision operator
$$Q(\mu_3^{1/2}\varphi_{k,0,0},\mu_3^{1/2}\varphi_{l,0,0})=a_{k,l}\ \mu_3^{1/2}\varphi_{k+l,0,0}, \quad a_{k,l} \in \rr, \quad k,l \geq 0,$$
which account for the triangularity of the system of differential equations (\ref{evv4}).

\bigskip

\begin{lemma}\label{proposition111}
Let $\varphi_{n,l,m} \in L^2(\rr^3)$ be the functions defined in \emph{(\ref{eig11})}.
The following algebraic identities hold:
\begin{itemize}
\item[$(i)$] $\mu_3^{-1/2}Q(\mu_3^{1/2}\varphi_{0,0,0},\mu_3^{1/2}\varphi_{0,0,0})=0$
\item[$(ii)$] $\displaystyle \mu_3^{-1/2}Q(\mu_3^{1/2}\varphi_{0,0,0},\mu_3^{1/2}\varphi_{m,0,0})=-\Big(\int_{-\frac\pi{4}}^{\frac\pi{4}}
\beta(\theta) \big(1-(\cos{\theta})^{2m}\big)d\theta\Big)\varphi_{m,0,0},$\\ $m \geq 1$
\item[$(iii)$] $\displaystyle \mu_3^{-1/2}Q(\mu_3^{1/2}\varphi_{n,0,0},\mu_3^{1/2}\varphi_{m,0,0})$\\
$\displaystyle=\sqrt{\frac{2n+2m+1}{(2n+1)(2m+1)}}\sqrt{C_{2n+2m}^{2n}}\Big(\int_{-\frac\pi{4}}^{\frac\pi{4}}
\beta(\theta)(\sin \theta)^{2n} (\cos \theta)^{2m}d\theta\Big)\varphi_{n+m,0,0},$\\
$n \geq 1,\ m\geq 0$
\end{itemize}
\end{lemma}

\bigskip

\noindent
In particular, this result allows to recover the diagonalization of the linearized radially symmetric Boltzmann operator
\begin{equation}\label{evvv3}
\mathscr{L}_1\varphi_{0,0,0}=0, \quad  \mathscr{L}_1\varphi_{n,0,0}=\Big(\int_{-\frac\pi{4}}^{\frac\pi{4}}\beta(\theta)(1-(\cos \theta)^{2n})d\theta\Big)\varphi_{n,0,0}, \quad n \geq 1,
\end{equation}
\begin{equation}\label{evvv4}
\mathscr{L}_2\varphi_{0,0,0}=0,   \quad \mathscr{L}_2\varphi_{n,0,0}=-\Big(\int_{-\frac\pi{4}}^{\frac\pi{4}}\beta(\theta)(\sin \theta)^{2n}d\theta\Big)\varphi_{n,0,0}, \quad n \geq 1,
\end{equation}
since $\varphi_{0,0,0}=\mu_3^{1/2}$.

\bigskip

\begin{proof}
We deduce from the Bobylev formula (Lemma~\ref{3.lem.kn223}) that for $m,n \geq 0$, $\xi \in \rr^3$,
\begin{multline*}
\mathcal{F}\big(Q(\mu_3^{1/2}\varphi_{n,0,0},\mu_3^{1/2}\varphi_{m,0,0})\big)(\xi) \\ = \int_{\val \theta\le \frac{\pi}{4}}\beta(\theta)\big[\widehat{\mu_3^{1/2}\varphi_{n,0,0}}(\xi \sin \theta)\widehat{\mu_3^{1/2}\varphi_{m,0,0}}(\xi\cos \theta)-\widehat{\mu_3^{1/2}\varphi_{n,0,0}}(0)\widehat{\mu_3^{1/2}\varphi_{m,0,0}}(\xi)
\big]d\theta.
\end{multline*}
Then, it follows from Lemma~\ref{eigbolt} that
\begin{multline*}
\mathcal{F}\big(Q(\mu_3^{1/2}\varphi_{n,0,0},\mu_3^{1/2}\varphi_{m,0,0})\big)(\xi) = \frac{(-1)^{m+n}}{\sqrt{(2n+1)(2m+1)}} \\ \times \int_{\val \theta\le \frac{\pi}{4}}\beta(\theta)\big[\widehat{\mu_1^{1/2}\psi_{2n}}(|\xi \sin \theta|)\widehat{\mu_1^{1/2}\psi_{2m}}(|\xi \cos \theta|)-(-1)^m\sqrt{2m+1}\widehat{\mu_1^{1/2}\psi_{2n}}(0)\widehat{\mu_3^{1/2}\varphi_{m,0,0}}(\xi)
\big]d\theta.
\end{multline*}
We deduce from (\ref{rumi1}) and Lemma~\ref{eigbolt} that
\begin{align*}\label{fh1}
& \ \widehat{\mu_1^{1/2}\psi_{2n}}(|\xi \sin \theta|) \widehat{\mu_1^{1/2}\psi_{2m}}(|\xi \cos \theta|)=  \frac{(-1)^{n+m}}{\sqrt{(2n)!(2m)!}}|\xi|^{2n+2m}(\sin \theta)^{2n} (\cos \theta)^{2m} e^{-\frac{|\xi|^2}{2}}\\ \notag
=  & \ \sqrt{C_{2n+2m}^{2n}}(\sin \theta)^{2n} (\cos \theta)^{2m}\widehat{\mu_1^{1/2}\psi_{2n+2m}}(|\xi|) \\ \notag
= & \ (-1)^{n+m} \sqrt{2n+2m+1} \sqrt{C_{2n+2m}^{2n}}(\sin \theta)^{2n} (\cos \theta)^{2m}\widehat{\mu_3^{1/2}\varphi_{n+m,0,0}}(\xi).
\end{align*}
Since
\begin{equation}\label{fh2}
\widehat{\mu_1^{1/2}\psi_{2n}}(0)=(\mu_1^{1/2},\psi_{2n})_{L^2}=(\psi_0,\psi_{2n})_{L^2}=\delta_{n,0},
\end{equation}
we obtain that
$$\mathcal{F}\big(Q(\mu_3^{1/2}\varphi_{n,0,0},\mu_3^{1/2}\varphi_{m,0,0})\big)(\xi)
=\alpha_{2n,2m}\sqrt{\frac{(2n+2m+1)}{(2n+1)(2m+1)}}\widehat{\mu_3^{1/2}\varphi_{n+m,0,0}}(\xi),$$
where the coefficients $\alpha_{2n,2m}$ are defined in (\ref{ar1}) and (\ref{ar2}).
This ends the proof of Lemma~\ref{proposition111}.
\end{proof}

\subsection{Trilinear estimates for the radially symmetric Boltzmann operator}
The key element and the main novelty of the present work is based on the derivation of sharp trilinear estimates for the terms
$$(\mu_3^{-1/2}Q(\mu_3^{1/2}f,\mu_3^{1/2}g),h)_{L^2}, \quad f,g,h \in \mathscr{S}_r(\rr^3).$$
The following lemma is instrumental in the proof of these trilinear estimates:

\bigskip

\begin{lemma}\label{lemkk000}
There exists a positive constant $C>0$ such that for all $n \geq 1$, $m \geq 0$,
$$\sqrt{\frac{2n+2m+1}{(2n+1)(2m+1)}}\sqrt{C_{2n+2m}^{2n}}\int_{-\frac{\pi}{4}}^{\frac{\pi}{4}}\beta(\theta)(\sin \theta)^{2n}(\cos \theta)^{2m}d\theta \leq  \frac{C}{n^{\frac{3}{4}}}\tilde{\mu}_{n,m},$$
where
$$\tilde{\mu}_{n,m}=\Big(1+\frac{m}{n}\Big)^{s}.$$
\end{lemma}

\bigskip

\begin{proof}
Lemma~\ref{lemkk000} is a direct consequence of the following estimates:
\begin{itemize}
\item[$(i)$] $\displaystyle \sqrt{C_{2n+2m}^{2n}}\int_{-\frac{\pi}{4}}^{\frac{\pi}{4}}\beta(\theta)(\sin \theta)^{2n}(\cos \theta)^{2m}d\theta \lesssim  \frac{1}{n^{\frac{3}{4}}} \Big(1+\frac{m}{n}\Big)^{s},$
when $m \gg 1$, $n \gg 1$
\item[$(ii)$] $\displaystyle \sqrt{C_{2n+2m}^{2n}}\int_{-\frac{\pi}{4}}^{\frac{\pi}{4}}\beta(\theta)(\sin \theta)^{2n}(\cos \theta)^{2m}d\theta \lesssim  m^{s},$
when $m \gg 1$, $1 \leq n \leq n_0$
\item[$(iii)$] $\displaystyle \sqrt{C_{2n+2m}^{2n}}\int_{-\frac{\pi}{4}}^{\frac{\pi}{4}}\beta(\theta)(\sin \theta)^{2n}(\cos \theta)^{2m}d\theta \lesssim  \frac{1}{n},$ when $0 \leq m \leq m_0$, $n \gg 1$
\end{itemize}
In order to establish these estimates, we begin by noticing that
$$\Lambda_{n,m}=\int_{-\frac{\pi}{4}}^{\frac{\pi}{4}}\beta(\theta)(\sin \theta)^{2n}(\cos \theta)^{2m}d\theta \approx \int_{0}^{\frac{\pi}{4}}(\sin \theta)^{2n-1-2s}(\cos \theta)^{2m+1}d\theta,$$
when the cross section $\beta$ satisfies 
$$\exists c>0, \quad  \frac{c}{|\theta|^{1+2s}} \leq \beta(\theta) \leq \frac{1}{c|\theta|^{1+2s}},$$
on the set $[-\frac{\pi}{4},\frac{\pi}{4}] \setminus \{0\}$.
By using the substitution rule with $t=\sin^2 \theta$, we obtain that
$$\int_{0}^{\frac{\pi}{4}}(\sin \theta)^{2n-1-2s}(\cos \theta)^{2m+1}d\theta=\frac{1}{2}\int_0^{\frac{1}{2}}t^{n-1-s}(1-t)^mdt.$$
This implies that
$$ \Lambda_{n,m} \approx \int_0^{\frac{1}{2}}t^{n-1-s}(1-t)^mdt.$$
By recalling the identity satisfied by the beta function
$$B(x,y)=\int_0^1t^{x-1}(1-t)^{y-1}dt=\frac{\Gamma(x)\Gamma(y)}{\Gamma(x+y)}, \quad \textrm{Re }x>0, \ \textrm{Re }y>0,$$
we obtain that
$$\int_0^{\frac{1}{2}}t^{n-1-s}(1-t)^mdt \leq B(n-s,m+1)=\frac{\Gamma(n-s)\Gamma(m+1)}{\Gamma(m+n+1-s)}.$$
By using the Stirling equivalent
$$\Gamma(x+1) \sim_{x \to +\infty} \sqrt{2\pi x}\Big(\frac{x}{e}\Big)^x, \quad \Gamma(n+1)=n!,$$
it follows that
$$\sqrt{C_{2n+2m}^{2n}}\Lambda_{n,m} \lesssim \sqrt{\frac{(2n+2m)!}{(2n)!(2m)!}}\frac{\Gamma(n-s)\Gamma(m+1)}{\Gamma(m+n+1-s)}.$$
This implies that
\begin{multline*}
\sqrt{C_{2n+2m}^{2n}}\Lambda_{n,m} \lesssim \Big(\frac{n+m}{n m}\Big)^{\frac{1}{4}}\Big(\frac{2n+2m}{e}\Big)^{n+m}\Big(\frac{e}{2n}\Big)^{n}\Big(\frac{e}{2m}\Big)^{m}\\
\times \sqrt{\frac{n m}{n+m}}\Big(\frac{n-s-1}{e}\Big)^{n-s-1}\Big(\frac{m}{e}\Big)^{m}\Big(\frac{e}{m+n-s}\Big)^{m+n-s},
\end{multline*}
when $m \gg 1$, $n \gg 1$. It follows that
\begin{align*}
\sqrt{C_{2n+2m}^{2n}}\Lambda_{n,m} \lesssim & \  \Big(\frac{n m}{n+m}\Big)^{\frac{1}{4}}\frac{(n+m)^{n+m}(n-s-1)^{n-s-1}}{(m+n-s)^{m+n-s}n^n} \\
 \lesssim & \ \Big(\frac{n m}{n+m}\Big)^{\frac{1}{4}}\frac{(m+n)^{s}}{(n-s-1)^{s+1}}\Big(1-\frac{s+1}{n}\Big)^{n}\Big(1+\frac{s}{m+n-s}\Big)^{m+n-s}\\
  \lesssim & \ m^{\frac{1}{4}}\frac{(m+n)^{s-\frac{1}{4}}}{n^{s+\frac{3}{4}}}  \lesssim  \frac{1}{n^{\frac{3}{4}}} \Big(1+\frac{m}{n}\Big)^{s},
\end{align*}
when $m \gg 1$, $n \gg 1$, since
\begin{equation}\label{ekkk0}
\forall n \geq 1, \forall x >-n, \quad \Big(1+\frac{x}{n}\Big)^n \leq e^x.
\end{equation}
When $m \gg 1$, $1 \leq n \leq n_0$, it follows that
$$\sqrt{C_{2n+2m}^{2n}}\Lambda_{n,m} \lesssim \sqrt{\frac{(2n+2m)!}{(2m)!}}\frac{\Gamma(m+1)}{\Gamma(m+n+1-s)}.$$
This implies that
\begin{multline*}
\sqrt{C_{2n+2m}^{2n}}\Lambda_{n,m} \lesssim \Big(\frac{n+m}{ m}\Big)^{\frac{1}{4}}\Big(\frac{2n+2m}{e}\Big)^{n+m}\Big(\frac{e}{2m}\Big)^{m}\\
\times \sqrt{\frac{m}{n+m-s}}\Big(\frac{m}{e}\Big)^{m}\Big(\frac{e}{m+n-s}\Big)^{m+n-s},
\end{multline*}
when $m \gg 1$, $1 \leq n \leq n_0$. We deduce from (\ref{ekkk0}) that
$$\sqrt{C_{2n+2m}^{2n}}\Lambda_{n,m} \lesssim   \frac{(n+m)^{n+m}}{(n+m-s)^{m+n-s}}=\Big(1+\frac{s}{n+m-s}\Big)^{m+n-s}(n+m)^{s} \lesssim m^s,$$
when $m \gg 1$, $1 \leq n \leq n_0$.
When $0 \leq m \leq m_0$, $n \gg 1$, it follows that
$$\sqrt{C_{2n+2m}^{2n}}\Lambda_{n,m} \lesssim \sqrt{\frac{(2n+2m)!}{(2n)!}}\frac{\Gamma(n-s)}{\Gamma(m+n+1-s)}.$$
This implies that
\begin{multline*}
\sqrt{C_{2n+2m}^{2n}}\Lambda_{n,m} \lesssim \Big(\frac{n+m}{n}\Big)^{\frac{1}{4}}\Big(\frac{2n+2m}{e}\Big)^{n+m}\Big(\frac{e}{2n}\Big)^{n}\\
\times \sqrt{\frac{n-s-1}{n+m-s}}\Big(\frac{n-s-1}{e}\Big)^{n-s-1}\Big(\frac{e}{m+n-s}\Big)^{m+n-s},
\end{multline*}
when $0 \leq m \leq m_0$, $n \gg 1$. We deduce from (\ref{ekkk0}) that
\begin{multline*}
\sqrt{C_{2n+2m}^{2n}}\Lambda_{n,m} \lesssim  \frac{(n+m)^{n+m}(n-s-1)^{n-s-1}}{n^n(m+n-s)^{m+n-s}} \\ \lesssim   \frac{(n+m)^{m}}{(n-s-1)^{m+1}}
\Big(1+\frac{m}{n}\Big)^{n}\Big(1-\frac{m+1}{m+n-s}\Big)^{m+n-s} \lesssim \frac{1}{n},
\end{multline*}
when $0 \leq m \leq m_0$, $n \gg 1$.
\end{proof}

\bigskip

\noindent
Lemma~\ref{lemkk000} allows to derive the following sharp trilinear estimates for the radially symmetric Boltzmann operator:

\bigskip

\begin{lemma}\label{proposition222}
There exists a positive constant $C>0$ such that for all $f,g,h \in \mathscr{S}_r(\rr^3)$,
$$|(\mu_3^{-1/2}Q(\mu_3^{1/2}f,\mu_3^{1/2}g),h)_{L^2}| \leq C \|f\|_{L^2}\|\mathcal{H}^\frac{s}{2}g\|_{L^2}\|\mathcal{H}^\frac{s}{2}h\|_{L^2},$$
and such that for all $f,g,h \in \mathscr{S}_r(\rr^3) \cap \mathcal{N}^{\perp}$, $n \geq 2$, $t \geq 0$,
\begin{multline*}
|(\mu_3^{-1/2}Q(\mu_3^{1/2}f,\mu_3^{1/2}g),e^{t\mathscr{L}}{\bf S}_nh)_{L^2}|
\leq C  \|e^{\frac{t}{2}\mathscr{L}}{\bf S}_{n-2}f\|_{L^2}\|e^{\frac{t}{2}\mathscr{L}}\mathcal{H}^\frac{s}{2}{\bf S}_{n-2}g\|_{L^2} \|e^{\frac{t}{2}\mathscr{L}}\mathcal{H}^\frac{s}{2}{\bf S}_{n}h\|_{L^2},
\end{multline*}
where $\mathscr{L}$ is the linearized non-cutoff Boltzmann operator, $\mathcal H=-\Delta_v+\frac{|v|^2}{4}$ is the 3-dimensional harmonic oscillator and ${\bf S}_n$ is the orthogonal projector onto the $n+1$ lowest energy levels
$${\bf S}_nf=\sum_{k=0}^n(f,\varphi_{k,0,0})_{L^2}\varphi_{k,0,0}, \quad e^{\frac{t}{2}\mathscr{L}}{\bf S}_nf=\sum_{k=0}^ne^{\frac{1}{2}\lambda_{2k}t}(f,\varphi_{k,0,0})_{L^2}\varphi_{k,0,0}.$$
\end{lemma}

\bigskip

\begin{proof}
Let $f,g,h \in \mathscr{S}_r(\rr^3)$ be some radial Schwartz functions. According to (\ref{rumi10}), we may decompose these functions into the Hermite basis $(\varphi_{n,l,m})_{n,l \geq 0, |m| \leq l}$ as follows
$$f=\sum_{n=0}^{+\infty}f_n \varphi_{n,0,0}, \quad g=\sum_{n=0}^{+\infty}g_n \varphi_{n,0,0}, \quad h=\sum_{n=0}^{+\infty}h_n \varphi_{n,0,0}.$$
We deduce from Lemma~\ref{proposition111} that
\begin{multline*}
(\mu_3^{-1/2}Q(\mu_3^{1/2}f,\mu_3^{1/2}g),h)_{L^2}=-f_{0} \sum_{n=0}^{+\infty}\overline{h_n} g_n\Big(\int_{-\frac\pi{4}}^{\frac\pi{4}}
\beta(\theta)(1-(\cos{\theta})^{2n})d\theta\Big)\\ +\sum_{n=0}^{+\infty}\overline{h_n}\Big(\sum_{\substack{k+l=n \\ k \geq 1, l \geq 0}}f_{k}g_l\sqrt{\frac{2k+2l+1}{(2k+1)(2l+1)}}\sqrt{C_{2k+2l}^{2k}}\int_{-\frac\pi{4}}^{\frac\pi{4}}
\beta(\theta)(\sin \theta)^{2k} (\cos \theta)^{2l}d\theta\Big).
\end{multline*}
It follows from Lemma~\ref{lemkk000} and (\ref{ekk1}) that
\begin{multline*}
|(\mu_3^{-1/2}Q(\mu_3^{1/2}f,\mu_3^{1/2}g),h)_{L^2}| \lesssim |f_{0}| \sum_{n=0}^{+\infty}\Big(2n+\frac{3}{2}\Big)^s|h_n| |g_n|\\ + \sum_{n=0}^{+\infty}|h_n|\Big(\sum_{\substack{k+l=n \\ k \geq 1, l \geq 0}}
\frac{\tilde{\mu}_{k,l}}{k^{\frac{3}{4}}} |f_{k}||g_l|\Big).
\end{multline*}
We have
\begin{align*}
|f_{0}| \sum_{n=0}^{+\infty}\Big(2n+\frac{3}{2}\Big)^s|h_n| |g_n| \lesssim & \ \|f\|_{L^2} \Big(\sum_{n=0}^{+\infty}\Big(2n+\frac{3}{2}\Big)^s|h_n|^2\Big)^{\frac{1}{2}}\Big(\sum_{n=0}^{+\infty}\Big(2n+\frac{3}{2}\Big)^s|g_n|^2\Big)^{\frac{1}{2}} \\
\lesssim & \ \|f\|_{L^2}\|\mathcal{H}^\frac{s}{2}g\|_{L^2}\|\mathcal{H}^\frac{s}{2}h\|_{L^2}.
\end{align*}
Furthermore, we notice that
\begin{align*}
 \sum_{n=0}^{+\infty}|h_n|\Big(\sum_{\substack{k+l=n \\ k \geq 1,l \geq 0}}\frac{\tilde{\mu}_{k,l}}{k^{\frac{3}{4}}} |f_{k}||g_l|\Big)=& \ \sum_{k \geq 1, l \geq 0}\frac{\tilde{\mu}_{k,l}}{k^{\frac{3}{4}}} |f_{k}||g_l||h_{k+l}|\\ =& \ \sum_{l=0}^{+\infty}\Big(2l+\frac{3}{2}\Big)^{\frac{s}{2}}|g_l| \Big(\sum_{k =1}^{+\infty}\frac{\tilde{\mu}_{k,l}}{k^{\frac{3}{4}}(2l+\frac{3}{2})^{\frac{s}{2}}} |f_{k}||h_{k+l}|\Big) .
\end{align*}
We deduce that
$$\sum_{n=0}^{+\infty}|h_n|\Big(\sum_{\substack{k+l=n \\ k \geq 1,l \geq 0}}\frac{\tilde{\mu}_{k,l}}{k^{\frac{3}{4}}} |f_{k}||g_l|\Big) \lesssim \|\mathcal{H}^\frac{s}{2}g\|_{L^2} \Big[\sum_{l=0}^{+\infty}\Big(\sum_{k =1}^{+\infty}\frac{\tilde{\mu}_{k,l}}{k^{\frac{3}{4}}(2l+\frac{3}{2})^{\frac{s}{2}}} |f_{k}||h_{k+l}|\Big)^2\Big]^{\frac{1}{2}}.$$
It follows that
\begin{align*}
 & \ \sum_{n=0}^{+\infty}|h_n|\Big(\sum_{\substack{k+l=n \\ k \geq 1,l \geq 0}}\frac{\tilde{\mu}_{k,l}}{k^{\frac{3}{4}}} |f_{k}||g_l|\Big)\\
 \lesssim & \ \|\mathcal{H}^\frac{s}{2}g\|_{L^2} \Big[\sum_{l=0}^{+\infty}\Big(\sum_{k =1}^{+\infty}|f_{k}|^2\Big)\Big(\sum_{k =1}^{+\infty}\frac{\tilde{\mu}_{k,l}^2}{k^{\frac{3}{2}}(2l+\frac{3}{2})^{s}}|h_{k+l}|^2\Big)\Big]^{\frac{1}{2}}\\
 \lesssim & \   \|f\|_{L^2}\|\mathcal{H}^\frac{s}{2}g\|_{L^2} \Big(\sum_{l=0}^{+\infty}\sum_{k =1}^{+\infty}\frac{\tilde{\mu}_{k,l}^2}{k^{\frac{3}{2}}(2l+\frac{3}{2})^{s}}|h_{k+l}|^2\Big)^{\frac{1}{2}}.
\end{align*}
We may write that
$$\Big(\sum_{l=0}^{+\infty}\sum_{k =1}^{+\infty}\frac{\tilde{\mu}_{k,l}^2}{k^{\frac{3}{2}}(2l+\frac{3}{2})^{s}}|h_{k+l}|^2\Big)^{\frac{1}{2}}=\Big[\sum_{n=0}^{+\infty}|h_{n}|^2\Big(\sum_{\substack{k+l =n\\ k\geq 1,l \geq 0}}\frac{\tilde{\mu}_{k,l}^2}{k^{\frac{3}{2}}(2l+\frac{3}{2})^{s}}\Big)\Big]^{\frac{1}{2}}.$$
On the other hand, we have
\begin{multline}\label{ekk233}
 \sum_{\substack{k+l =n\\ k\geq 1,l \geq 0}}\frac{\tilde{\mu}_{k,l}^2}{k^{\frac{3}{2}}(2l+\frac{3}{2})^{s}} \lesssim \sum_{\substack{k+l =n\\ k\geq 1,l \geq 0\\ k \geq l}}\frac{1}{k^{\frac{3}{2}}(2l+\frac{3}{2})^{s}}+\sum_{\substack{k+l =n\\ k\geq 1,l \geq 0 \\ k \leq l}}\frac{(2l+\frac{3}{2})^{s}}{k^{\frac{3}{2}}}
\\ \lesssim  \sum_{k= 1}^{n}\frac{1}{k^{\frac{3}{2}}}+\Big(2n+\frac{3}{2}\Big)^{s}\sum_{k= 1}^{n}\frac{1}{k^{\frac{3}{2}}}
\lesssim   \Big[\Big(2n+\frac{3}{2}\Big)^{s}+1\Big]\sum_{k= 1}^{+\infty}\frac{1}{k^{\frac{3}{2}}}\lesssim \Big(2n+\frac{3}{2}\Big)^{s},
\end{multline}
since
$$\tilde{\mu}_{k,l} \lesssim 1 \textrm{ when } k \geq l,\ k \geq 1, \quad \tilde{\mu}_{k,l} \lesssim \Big(2l+\frac{3}{2}\Big)^{s} \textrm{ when } 1 \leq k \leq l.$$
This implies that
$$\Big(\sum_{l=0}^{+\infty}\sum_{k =1}^{+\infty}\frac{\tilde{\mu}_{k,l}^2}{k^{\frac{3}{2}}(2l+\frac{3}{2})^{s}}|h_{k+l}|^2\Big)^{\frac{1}{2}} \lesssim \|\mathcal{H}^\frac{s}{2}h\|_{L^2}.$$
We conclude that
$$|(\mu_3^{-1/2}Q(\mu_3^{1/2}f,\mu_3^{1/2}g),h)_{L^2}| \lesssim \|f\|_{L^2}\|\mathcal{H}^\frac{s}{2}g\|_{L^2}\|\mathcal{H}^\frac{s}{2}h\|_{L^2}.$$
On the other hand, we may write that for all $f,g,h \in \mathscr{S}_r(\rr^3)$, $n \geq 2$,
\begin{align*}
& \ (\mu_3^{-1/2}Q(\mu_3^{1/2}f,\mu_3^{1/2}g),e^{t\mathscr{L}}{\bf S}_nh)_{L^2}=\sum_{j=1}^{n}e^{\lambda_{2j}t}\overline{h_j}f_j g_0\int_{-\frac\pi{4}}^{\frac\pi{4}}
\beta(\theta)(\sin \theta)^{2j}d\theta\\
+& \ \sum_{j=2}^{n}e^{\lambda_{2j}t}\overline{h_j}\Big(\sum_{\substack{k+l=j \\ k \geq 1, l \geq 1}}f_{k}g_l\sqrt{\frac{2k+2l+1}{(2k+1)(2l+1)}}\sqrt{C_{2k+2l}^{2k}}\int_{-\frac\pi{4}}^{\frac\pi{4}}
\beta(\theta)(\sin \theta)^{2k} (\cos \theta)^{2l}d\theta\Big)\\
- & \ f_{0} \sum_{j=0}^{n}e^{\lambda_{2j}t}\overline{h_j} g_j\int_{-\frac\pi{4}}^{\frac\pi{4}}
\beta(\theta)(1-(\cos{\theta})^{2j})d\theta.
\end{align*}
We assume that  $f,g,h \in \mathscr{S}_r(\rr^3) \cap \mathcal{N}^{\perp}$, i.e.
$$f_0=f_1=g_0=g_1=h_0=h_1=0.$$
This implies that
$$|(\mu_3^{-1/2}Q(\mu_3^{1/2}f,\mu_3^{1/2}g),e^{t\mathscr{L}}{\bf S}_nh)_{L^2}| \lesssim \sum_{j=4}^{n}e^{\lambda_{2j}t}|h_j|\Big(\sum_{\substack{k+l=j \\ k \geq 2, l \geq 2}}
\frac{\tilde{\mu}_{k,l}}{k^{\frac{3}{4}}} |f_{k}||g_l|\Big).$$
We notice that
\begin{multline*}
 \sum_{j=4}^{n}e^{\lambda_{2j}t}|h_j|\Big(\sum_{\substack{k+l=j \\ k \geq 2,l \geq 2}}\frac{\tilde{\mu}_{k,l}}{k^{\frac{3}{4}}} |f_{k}||g_l|\Big)= \sum_{\substack{k \geq 2, l \geq 2 \\ 4 \leq k+l \leq n}}\frac{\tilde{\mu}_{k,l}}{k^{\frac{3}{4}}}e^{\lambda_{2k+2l}t} |f_{k}||g_l||h_{k+l}|\\ = \sum_{l=2}^{n-2}\Big(2l+\frac{3}{2}\Big)^{\frac{s}{2}}e^{\frac{\lambda_{2l}}{2}t}|g_l| \Big(\sum_{\substack{k \geq 2 \\ 4 \leq k+l \leq n}}  \frac{e^{-\frac{t}{2}(\lambda_{2k}+\lambda_{2l}-\lambda_{2k+2l})}\tilde{\mu}_{k,l}}{k^{\frac{3}{4}}(2l+\frac{3}{2})^{\frac{s}{2}}} e^{\frac{\lambda_{2k}}{2}t}|f_{k}|e^{\frac{\lambda_{2k+2l}}{2}t}|h_{k+l}|\Big) .
\end{multline*}
Since by Lemma~\ref{lemk10}, we have
$$\forall t \geq 0, \forall k,l \geq 2, \quad  e^{-\frac{t}{2}(\lambda_{2k}+\lambda_{2l}-\lambda_{2k+2l})} \leq 1,$$
we obtain that
\begin{multline*}
|(\mu_3^{-1/2}Q(\mu_3^{1/2}f,\mu_3^{1/2}g),e^{t\mathscr{L}}{\bf S}_nh)_{L^2}|  \\
\leq \sum_{l=2}^{n-2}\Big(2l+\frac{3}{2}\Big)^{\frac{s}{2}}e^{\frac{\lambda_{2l}}{2}t}|g_l| \Big(\sum_{\substack{k \geq 2 \\ 4 \leq k+l \leq n}} \frac{\tilde{\mu}_{k,l}}{k^{\frac{3}{4}}(2l+\frac{3}{2})^{\frac{s}{2}}} e^{\frac{\lambda_{2k}}{2}t}|f_{k}|e^{\frac{\lambda_{2k+2l}}{2}t}|h_{k+l}|\Big)
 \\ \leq \|e^{\frac{t}{2}\mathscr{L}}\mathcal{H}^\frac{s}{2}\textbf{S}_{n-2}g\|_{L^2} \Big[\sum_{l=2}^{n-2}\Big(\sum_{\substack{k \geq 2 \\ 4 \leq k+l \leq n}}\frac{\tilde{\mu}_{k,l}}{k^{\frac{3}{4}}(2l+\frac{3}{2})^{\frac{s}{2}}} e^{\frac{\lambda_{2k}}{2}t}|f_{k}|e^{\frac{\lambda_{2k+2l}}{2}t}|h_{k+l}|\Big)^2\Big]^{\frac{1}{2}}.
\end{multline*}
It follows that
\begin{align*}
 & \ |(\mu_3^{-1/2}Q(\mu_3^{1/2}f,\mu_3^{1/2}g),e^{t\mathscr{L}}{\bf S}_nh)_{L^2}| \\
 \leq & \ \|e^{\frac{t}{2}\mathscr{L}}\mathcal{H}^\frac{s}{2}{\bf S}_{n-2}g\|_{L^2} \Big[\sum_{l=2}^{n-2}\Big(\sum_{k=2}^{n-2}e^{\lambda_{2k}t}|f_{k}|^2\Big)\Big(\sum_{\substack{k \geq 2 \\ 4 \leq k+l \leq n}}\frac{\tilde{\mu}_{k,l}^2}{k^{\frac{3}{2}}(2l+\frac{3}{2})^{s}}e^{\lambda_{2k+2l}t}|h_{k+l}|^2\Big)\Big]^{\frac{1}{2}}\\
 \leq & \   \|e^{\frac{t}{2}\mathscr{L}}{\bf S}_{n-2}f\|_{L^2}\|e^{\frac{t}{2}\mathscr{L}}\mathcal{H}^\frac{s}{2}{\bf S}_{n-2}g\|_{L^2} \Big(\sum_{\substack{k \geq 2, l \geq 2 \\ 4 \leq k+l \leq n}}\frac{\tilde{\mu}_{k,l}^2}{k^{\frac{3}{2}}(2l+\frac{3}{2})^{s}}e^{\lambda_{2k+2l}t}|h_{k+l}|^2\Big)^{\frac{1}{2}}.
\end{align*}
By arguing as before, we deduce from (\ref{ekk233}) that
\begin{multline*}
\Big(\sum_{\substack{k \geq 2, l \geq 2 \\ 4 \leq k+l \leq n}}\frac{\tilde{\mu}_{k,l}^2}{k^{\frac{3}{2}}(2l+\frac{3}{2})^{s}}e^{\lambda_{2k+2l}t}|h_{k+l}|^2\Big)^{\frac{1}{2}}=\Big[\sum_{j=4}^{n}e^{\lambda_{2j}t}|h_{j}|^2\Big(\sum_{\substack{k+l =j\\ k\geq 2,l \geq 2}}\frac{\tilde{\mu}_{k,l}^2}{k^{\frac{3}{2}}(2l+\frac{3}{2})^{s}}\Big)\Big]^{\frac{1}{2}} \\
\lesssim \Big[\sum_{j=4}^{n}\Big(2j+\frac{3}{2}\Big)^se^{\lambda_{2j}t} |h_{j}|^2\Big]^{\frac{1}{2}} \lesssim \|e^{\frac{t}{2}\mathscr{L}}\mathcal{H}^\frac{s}{2}{\bf S}_{n}h\|_{L^2}.
\end{multline*}
This implies that
\begin{multline*}
|(\mu_3^{-1/2}Q(\mu_3^{1/2}f,\mu_3^{1/2}g),e^{t\mathscr{L}}{\bf S}_nh)_{L^2}|
 \lesssim  \|e^{\frac{t}{2}\mathscr{L}}{\bf S}_{n-2}f\|_{L^2}\|e^{\frac{t}{2}\mathscr{L}}\mathcal{H}^\frac{s}{2}{\bf S}_{n-2}g\|_{L^2} \|e^{\frac{t}{2}\mathcal{L}}\mathcal{H}^\frac{s}{2}{\bf S}_{n}h\|_{L^2}.
\end{multline*}
This ends the proof of Lemma~\ref{proposition222}.
\end{proof}

\subsection{Explicit solutions to the radially symmetric spatially homogeneous Boltzmann equation with small initial radial data}
\noindent
We solve explicitly the Cauchy problem associated to the non-cutoff radially symmetric spatially homogeneous Boltzmann equation with Maxwellian molecules for a small $L^2$-initial radial data
$$\begin{cases}
\partial_tg+\mathscr{L}g=\mu_3^{-1/2}Q(\sqrt{\mu_3}g,\sqrt{\mu_3}g),\\
g|_{t=0}=g_0.
\end{cases}$$
We search for a radial solution
\begin{equation}\label{rumi20}
g(t)=\sum_{n=0}^{+\infty}b_n(t)\varphi_{n,0,0}.
\end{equation}
It follows from Lemma~\ref{proposition111} that
\begin{multline*}
\mu_3^{-1/2}Q(\sqrt{\mu_3}g,\sqrt{\mu_3}g)=-\sum_{k=1}^{+\infty}b_0(t)b_k(t)\Big(\int_{-\frac\pi{4}}^{\frac\pi{4}}
\beta(\theta) \big(1-(\cos{\theta})^{2k}\big)d\theta\Big)\varphi_{k,0,0}\\
+\sum_{k \geq 1, l \geq 0}b_{k}(t)b_l(t)\sqrt{\frac{2k+2l+1}{(2k+1)(2l+1)}}\sqrt{C_{2k+2l}^{2k}}\Big(\int_{-\frac\pi{4}}^{\frac\pi{4}}
\beta(\theta)(\sin \theta)^{2k} (\cos \theta)^{2l}d\theta\Big)\varphi_{k+l,0,0}.
\end{multline*}
This implies that
\begin{equation}\label{dd1}
\mu_3^{-1/2}Q(\sqrt{\mu_3}g,\sqrt{\mu_3}g)=\sum_{n=1}^{+\infty}\tilde{\alpha}_n(t)\varphi_{n,0,0},
\end{equation}
where
\begin{multline}\label{dd2}
\tilde{\alpha}_n(t)=-b_0(t)b_n(t)\Big(\int_{-\frac\pi{4}}^{\frac\pi{4}}\beta(\theta) \big(1-(\cos{\theta})^{2n}\big)d\theta\Big)\\
+\sum_{\substack{k+l=n \\ k \geq 1, l \geq 0}}b_{k}(t)b_l(t)\sqrt{\frac{2k+2l+1}{(2k+1)(2l+1)}}\sqrt{C_{2k+2l}^{2k}}\Big(\int_{-\frac\pi{4}}^{\frac\pi{4}}
\beta(\theta)(\sin \theta)^{2k} (\cos \theta)^{2l}d\theta\Big).
\end{multline}
By using that
\begin{equation}\label{hk1}
\mathscr{L}g=\sum_{n=2}^{+\infty}\lambda_{2n} b_n(t)\varphi_{n,0,0},
\end{equation}
we notice that the radially symmetric spatially homogeneous Boltzmann equation is satisfied if and only if the functions $(b_n)_{n\geq 0}$ satisfy  the infinite system of differential equations (\ref{evv3}),
$$\left\{
\begin{array}{c}
    \partial_t b_0(t) = 0 \\
 \displaystyle   \forall n \geq 1, \ \partial_t b_n(t)+\lambda_{2n} b_n(t) =   -b_0(t)b_n(t)\Big(\int_{-\frac\pi{4}}^{\frac\pi{4}}\beta(\theta) \big(1-(\cos{\theta})^{2n}\big)d\theta\Big)\\
 \displaystyle +\sum_{\substack{k+l=n \\ k \geq 1, l \geq 0}}b_{k}(t)b_l(t)\sqrt{\frac{2k+2l+1}{(2k+1)(2l+1)}}\sqrt{C_{2k+2l}^{2k}}\Big(\int_{-\frac\pi{4}}^{\frac\pi{4}}
\beta(\theta)(\sin \theta)^{2k} (\cos \theta)^{2l}d\theta\Big).
\end{array}
\right.$$
When $b_0(0)=b_1(0)=0$, this system reduces to
$$\left\{
\begin{array}{c}
\displaystyle \forall t \geq 0, \quad b_0(t)=b_1(t)=0, \\
\displaystyle \forall t \geq 0, \quad \partial_t b_2(t)+\lambda_4 b_2(t)=0,\\
\displaystyle \forall n \geq 3, \forall t \geq 0, \quad \partial_t b_n(t)+\lambda_{2n} b_n(t)\\
\displaystyle  \qquad =\sum_{\substack{k+l=n \\ k \geq 1, l \geq 1}}b_{k}(t)b_l(t)\sqrt{\frac{2k+2l+1}{(2k+1)(2l+1)}}\sqrt{C_{2k+2l}^{2k}}\Big(\int_{-\frac\pi{4}}^{\frac\pi{4}}
\beta(\theta)(\sin \theta)^{2k} (\cos \theta)^{2l}d\theta\Big).
\end{array}
\right.$$
As mentioned previously, this system of differential equations is triangular since the $(n+1)^{\textrm{th}}$ equation is a linear differential equation for the function $b_n$ with a right-hand-side involving only the functions $(b_k)_{1 \leq k \leq n-1}$. This system may therefore be explicitly solved while solving a sequence of linear differential equations.

Let $g_0 \in L^2(\rr)$ be a radial fluctuation belonging to the orthogonal of the space of collisional invariants $\mathcal{N}^{\perp},$
i.e.,
$$(g_0,\varphi_{0,0,0})_{L^2}=(g_0,\varphi_{1,0,0})_{L^2}=0,$$
and $(b_n(t))_{n \geq 0}$ be the solutions of the system (\ref{evv4}) with initial conditions
\begin{equation}\label{bkk4yyy}
\forall n \geq 0, \quad b_n(t)|_{t=0}=(g_0,\varphi_{n,0,0})_{L^2}.
\end{equation}
According to the theory of linear differential equations, these functions exist and are uniquely defined on $\rr$.
We consider the function
\begin{equation}\label{bkk3yyy}
g(t)=\sum_{n=2}^{+\infty}b_n(t)\varphi_{n,0,0}
\end{equation}
and
$$(\textbf{S}_Ng)(t)=\sum_{n=2}^{N}b_n(t)\varphi_{n,0,0}, \quad N \geq 2,$$
its orthogonal projection onto the $N+1$ lowest energy levels.
Thanks to the definition of the functions $(b_n)_{n \geq 0}$, we deduce from (\ref{dd1}) and (\ref{dd2}) that for all $N \geq 2$, $t \geq 0$,
\begin{multline*}
\big(\partial_t(\textbf{S}_Ng)(t),e^{t\mathscr{L}}(\textbf{S}_Ng)(t)\big)_{L^2}+\big(\mathscr{L}(\textbf{S}_Ng)(t),e^{t\mathscr{L}}(\textbf{S}_Ng)(t)\big)_{L^2}\\ =\big(\mu_3^{-1/2}Q\big(\sqrt{\mu_3}(\textbf{S}_Ng)(t),\sqrt{\mu_3}(\textbf{S}_Ng)(t)\big),e^{t\mathcal{L}}(\textbf{S}_Ng)(t)\big)_{L^2}.
\end{multline*}
Since the linearized Boltzmann operator is selfadjoint, we notice that
\begin{align*}
& \ \frac{1}{2}\frac{d}{dt}\|e^{\frac{t}{2}\mathscr{L}}(\textbf{S}_Ng)(t)\|_{L^2}^2+\frac{1}{2}\big(\mathscr{L}e^{\frac{t}{2}\mathscr{L}}(\textbf{S}_Ng)(t),e^{\frac{t}{2}\mathscr{L}}(\textbf{S}_Ng)(t)\big)_{L^2}\\
= & \ \textrm{Re}\big(\partial_t(\textbf{S}_Ng)(t),e^{t\mathscr{L}}(\textbf{S}_Ng)(t)\big)_{L^2}+\textrm{Re}\big(\mathscr{L}e^{\frac{t}{2}\mathscr{L}}(\textbf{S}_Ng)(t),e^{\frac{t}{2}\mathscr{L}}(\textbf{S}_Ng)(t)\big)_{L^2}\\
= & \ \textrm{Re}\big(\partial_t(\textbf{S}_Ng)(t),e^{t\mathscr{L}}(\textbf{S}_Ng)(t)\big)_{L^2}+\textrm{Re}\big(\mathscr{L}(\textbf{S}_Ng)(t),e^{t\mathscr{L}}(\textbf{S}_Ng)(t)\big)_{L^2}.
\end{align*}
It follows that
\begin{multline*}
\frac{1}{2}\frac{d}{dt}\|e^{\frac{t}{2}\mathscr{L}}(\textbf{S}_Ng)(t)\|_{L^2}^2+\frac{1}{2}\big(\mathscr{L}e^{\frac{t}{2}\mathscr{L}}(\textbf{S}_Ng)(t),e^{\frac{t}{2}\mathscr{L}}(\textbf{S}_Ng)(t)\big)_{L^2} \\
= \textrm{Re}\big(\mu_3^{-1/2}Q\big(\sqrt{\mu_3}(\textbf{S}_Ng)(t),\sqrt{\mu_3}(\textbf{S}_Ng)(t)\big),e^{t\mathcal{L}}(\textbf{S}_Ng)(t)\big)_{L^2}.
\end{multline*}
We deduce from Lemma~\ref{proposition222} that for all $N \geq 2$, $t \geq 0$,
\begin{multline*}
\frac{1}{2}\frac{d}{dt}\|e^{\frac{t}{2}\mathscr{L}}(\textbf{S}_{N}g)(t)\|_{L^2}^2+\frac{1}{2}\big(\mathscr{L}e^{\frac{t}{2}\mathscr{L}}(\textbf{S}_Ng)(t),e^{\frac{t}{2}\mathscr{L}}(\textbf{S}_Ng)(t)\big)_{L^2} \\
\leq   C  \|e^{\frac{t}{2}\mathscr{L}}({\bf S}_{N-2}g)(t)\|_{L^2}\|e^{\frac{t}{2}\mathscr{L}}\mathcal{H}^\frac{s}{2}({\bf S}_{N}g)(t)\|_{L^2}^2.
\end{multline*}
According to (\ref{hk1}), this implies that
\begin{multline*}
\frac{1}{2}\frac{d}{dt}\|e^{\frac{t}{2}\mathscr{L}}(\textbf{S}_{N}g)(t)\|_{L^2}^2+\frac{1}{2}\sum_{n=2}^{N}\lambda_{2n}e^{\lambda_{2n}t}|b_n(t)|^2 \\
\leq   C  \|e^{\frac{t}{2}\mathscr{L}}({\bf S}_{N-2}g)(t)\|_{L^2}\sum_{n=2}^{N}\Big(2n+\frac{3}{2}\Big)^se^{\lambda_{2n}t}|b_n(t)|^2.
\end{multline*}
It follows from (\ref{ekk1}) that we may find a new positive constant $C>0$ such that for all $N \geq 2$, $t \geq 0$,
\begin{multline}\label{hk2}
\frac{d}{dt}\|e^{\frac{t}{2}\mathscr{L}}(\textbf{S}_{N}g)(t)\|_{L^2}^2+\sum_{n=2}^{N}\lambda_{2n}e^{\lambda_{2n}t}|b_n(t)|^2 \\
\leq   C  \|e^{\frac{t}{2}\mathscr{L}}({\bf S}_{N-2}g)(t)\|_{L^2}\sum_{n=2}^{N}\lambda_{2n}e^{\lambda_{2n}t}|b_n(t)|^2.
\end{multline}

\bigskip

\begin{lemma}\label{lemk1111}
There exists a positive constant $\eps_0>0$ such that for all $0<\eps \leq \eps_0$,
$$g_0 \in \mathcal{N}^{\perp}, \ \|g_0\|_{L^2} \leq \eps \Rightarrow \forall N \geq 0, \forall t \geq 0, \quad \|e^{\frac{t}{2}\mathscr{L}}(\textbf{S}_{N}g)(t)\|_{L^2} \leq \eps.$$
\end{lemma}

\bigskip

\begin{proof}
When $0<\eps \leq \eps_0<1$ and $\|g_0\|_{L^2} \leq \eps$, we deduce from (\ref{kek5}) and (\ref{kek6}) that for all $t \geq 0$,
$$\|e^{\frac{t}{2}\mathscr{L}}(\textbf{S}_{0}g)(t)\|_{L^2}^2=|b_0(t)|^2=0, \quad \|e^{\frac{t}{2}\mathscr{L}}(\textbf{S}_{1}g)(t)\|_{L^2}^2=|b_0(t)|^2+|b_1(t)|^2=0,$$
\begin{multline*}
\|e^{\frac{t}{2}\mathscr{L}}(\textbf{S}_{2}g)(t)\|_{L^2}^2=|b_0(t)|^2+|b_1(t)|^2+e^{\lambda_4t}|b_2(t)|^2=|b_2(0)|^2e^{-\lambda_4 t} \leq |b_2(0)|^2
\leq \|g_0\|_{L^2}^2 \leq \eps^2.
\end{multline*}
We choose the positive constant $\eps_0>0$ such that
$$0<\eps_0 \leq \frac{1}{2C}.$$
With that choice, we notice that the condition
$$\forall t \geq 0, \quad \|e^{\frac{t}{2}\mathscr{L}}(\textbf{S}_{N-1}g)(t)\|_{L^2} \leq \eps,$$
implies first that for all $t \geq 0$,
\begin{multline*}
\|e^{\frac{t}{2}\mathscr{L}}(\textbf{S}_{N-2}g)(t)\|_{L^2}=\Big(\sum_{n=2}^{N-2}e^{\lambda_{2n}t}|b_n(t)|^2\Big)^{\frac{1}{2}} \\ \leq \Big(\sum_{n=2}^{N-1}e^{\lambda_{2n}t}|b_n(t)|^2\Big)^{\frac{1}{2}}= \|e^{\frac{t}{2}\mathscr{L}}(\textbf{S}_{N-1}g)(t)\|_{L^2} \leq \eps.
\end{multline*}
It therefore follows from (\ref{hk2}) that for all $t \geq 0$,
\begin{align}\label{ek7yyy}
& \ \frac{d}{dt}\|e^{\frac{t}{2}\mathscr{L}}(\textbf{S}_{N}g)(t)\|_{L^2}^2+\frac{1}{2}\sum_{n=2}^{N}\lambda_{2n}e^{t\lambda_{2n}}|b_n(t)|^2 \\
\leq & \ \frac{d}{dt}\|e^{\frac{t}{2}\mathscr{L}}(\textbf{S}_{N}g)(t)\|_{L^2}^2+(1-C\eps)\sum_{n=2}^{N}\lambda_{2n}e^{t\lambda_{2n}}|b_n(t)|^2
 \leq  0, \notag
\end{align}
thus
$$\forall t \geq 0, \quad \|e^{\frac{t}{2}\mathscr{L}}(\textbf{S}_{N}g)(t)\|_{L^2}  \leq \|(\textbf{S}_{N}g)(0)\|_{L^2} \leq \|g_0\|_{L^2} \leq \eps.$$
The proof of Lemma~\ref{lemk1111} follows by induction.
\end{proof}

\bigskip

\noindent
We may now finish the proof of Theorem~\ref{thnonlinear11}. The global well-posedness for small initial radial data is a straightforward consequence of Lemma~\ref{lemk1111} which ensures that there exists a positive constant $\eps_0>0$ such that the fluctuation defined in (\ref{rumi20}) satisfies 
$$\|g_0\|_{L^2} \leq \eps_0 \Rightarrow \forall t \geq 0, \ \|g(t)\|_{L^2} \leq \|e^{\frac{t}{2}\mathscr{L}}g(t)\|_{L^2} \leq  \eps_0 <+\infty.$$
The uniqueness is a direct consequence of the uniqueness of the solutions to the system of differential equations (\ref{evv4}) with initial conditions
$$\forall n \geq 0, \quad b_n(t)|_{t=0}=(g_0,\varphi_{n,0,0})_{L^2}.$$
Let $0<\delta<1$ be a positive constant. We deduce from Lemma~\ref{lemk1111} and (\ref{hk2}) that there exists a positive constant
$$0<\eps \leq \inf\Big(\frac{\delta}{C},\eps_0\Big),$$
such that, when $g_0 \in \mathcal{N}^{\perp}$, $\|g_0\|_{L^2} \leq \eps$, we have for all $N \geq 2$,
\begin{multline*}
\frac{d}{dt}\|e^{\frac{t}{2}\mathscr{L}}(\textbf{S}_{N}g)(t)\|_{L^2}^2+\lambda_{4}(1-\delta)\|e^{\frac{t}{2}\mathscr{L}}(\textbf{S}_{N}g)(t)\|_{L^2}^2
= \frac{d}{dt}\|e^{\frac{t}{2}\mathscr{L}}(\textbf{S}_{N}g)(t)\|_{L^2}^2\\ +\lambda_{4}(1-\delta)\sum_{n=2}^{N}e^{\lambda_{2n}t}|b_n(t)|^2
\leq \frac{d}{dt}\|e^{\frac{t}{2}\mathscr{L}}(\textbf{S}_{N}g)(t)\|_{L^2}^2+(1-\delta)\sum_{n=2}^{N}\lambda_{2n}e^{\lambda_{2n}t}|b_n(t)|^2 \leq 0,
\end{multline*}
since
$$\forall n \geq 2, \quad  \lambda_4 \leq \lambda_{2n}.$$
This implies that
$$\forall N \geq 2, \forall t \geq 0, \quad \|e^{\frac{t}{2}\mathscr{L}}(\textbf{S}_{N}g)(t)\|_{L^2} \leq e^{-\frac{\lambda_4}{2}(1-\delta) t}\|g_0\|_{L^2}.$$
We therefore deduce that
$$\forall t \geq 0, \quad \|e^{\frac{t}{2}\mathscr{L}}g(t)\|_{L^2} \leq e^{-\frac{\lambda_4}{2}(1-\delta) t}\|g_0\|_{L^2}.$$
This ends the proof of Theorem~\ref{thnonlinear11}.

\section{Appendix}\label{appendix}

\subsection{The Bobylev formula for the radially symmetric Boltzmann operator}\label{bobylev1}
The Bobylev formula is very useful identity providing an explicit formula for the Fourier transform of the Boltzmann operator \eqref{eq1}, see e.g.~\cite{al-1}. In the case of Maxwellian molecules, the Fourier transform of the Boltzmann operator whose cross section satisfies the assumption (\ref{sa1}),
$$Q(g, f)(v)=\int_{\rr^d}\int_{\SSS^{d-1}}b\Big(\frac{v-v_{*}}{|v-v_{*}|} \cdot \sigma\Big)\bigl(g'_* f'-g_{*}f\bigr)d\sigma dv_*,$$
is equal to
\begin{align*}
\mathcal{F}\big(Q(g, f)\big)(\xi)= & \ \int_{\rr^d}Q(g, f)(v)e^{-iv \cdot \xi}dv\\
= & \ \int_{\SSS^{d-1}} b\Big(\frac{\xi}{|\xi|} \cdot \sigma\Big)\bigl[\widehat{g}(\xi^-)\widehat{f}(\xi^+)-
\widehat{g}(0)\widehat{f}(\xi)\bigr]d\sigma,
\end{align*}
where $\xi^+=\frac{\xi+|\xi|\sigma}{2}$ and $\xi^-=\frac{\xi-|\xi|\sigma}{2}$.

\bigskip

\begin{lemma} \label{3.lem.kn223}
For $f,g\in \mathscr S_{r}(\R^3)$, we have
\begin{equation}\label{bob1}
\mathcal{F}\big(Q(g, f)\big)(\xi)=\int_{\val \theta\le \frac{\pi}{4}}\beta(\theta)\big[\hat g(\xi \sin \theta)\hat f(\xi\cos \theta)-\hat g(0)\hat f(\xi)
\big]d\theta,
\end{equation}
where $\beta$ is the function defined in \eqref{new001}.
\end{lemma}

\bigskip

\begin{proof}
Thanks to the Bobylev formula, we may write with $\nu=\frac{\xi}{\val \xi}$,
\begin{multline}\label{3.kkww2}
\mathcal{F}\big(Q(g, f)\big)(\xi)=\int_{(0,\pi)_{\theta}\times\mathbb S^{1}_{\omega}}
b(\cos \theta)\sin \theta \\
\times \Big[\hat g\left(\frac{\xi-\val \xi( \omega\sin \theta\oplus \nu\cos \theta)}2\right)\hat f\left(\frac{\xi+\val \xi( \omega\sin \theta\oplus \nu\cos \theta)}2\right)-\hat g(0)\hat f(\xi)\Big]d\theta d\omega.
\end{multline}
The cross section $b(\cos \theta)$ is supported where $0\le \theta\le \frac{\pi}{2}$ and we notice that
\begin{align*}
\xi-\val \xi( \omega\sin \theta\oplus \nu\cos \theta)=& \ \val\xi\bigl(-\omega\sin\theta\oplus\nu(1-\cos \theta)\bigr)\\
= & \ 2\val \xi\sin\Big(\frac\theta2\Big)\Big[-\omega\cos\Big(\frac\theta2\Big)\oplus\nu\sin\Big(\frac\theta2\Big)\Big],
\end{align*}
\begin{align*}
\xi+\val \xi( \omega\sin \theta\oplus \nu\cos \theta)= & \ \val\xi\bigl(\omega\sin\theta\oplus\nu(1+\cos \theta)\bigr)\\
= & \ 2\val \xi\cos\Big(\frac\theta2\Big)\Big[\omega\sin\Big(\frac\theta2\Big)\oplus\nu\cos\Big(\frac\theta2\Big)\Big],
\end{align*}
so that, since $\hat g, \hat f$ are radial functions,
$$\hat g\Big(\frac{\xi-\val \xi( \omega\sin \theta\oplus \nu\cos \theta)}2\Big)=\hat g\Big(\val \xi\sin\Big(\frac\theta2\Big) \nu\Big)=\hat g\Big(\xi\sin\Big(\frac\theta2\Big)\Big),$$
$$\hat f\Big(\frac{\xi+\val \xi( \omega\sin \theta\oplus \nu\cos \theta)}2\Big)=\hat f\Big(\val \xi\cos\Big(\frac\theta2\Big) \nu\Big)=\hat f\Big(\xi\cos\Big(\frac\theta2\Big)\Big),$$
yielding
{\small \begin{align*}
\mathcal{F}\big(Q(g, f)\big)&(\xi)=\val{\mathbb S^{1}}\int_{0}^{\frac{\pi}{2}}b(\cos \theta)\sin \theta\Big[\hat g\Big(\xi\sin\Big(\frac\theta2\Big)\Big)
\hat f\Big(\xi\cos\Big(\frac\theta2\Big)\Big)-\hat g(0)\hat f(\xi)\Big]d\theta\\
&=2\int_{0}^{\frac{\pi}{4}}\underbrace{\val{\mathbb S^{1}}b(\cos2\theta)\sin2\theta  }_{=\beta(\theta)\text{ from \eqref{new001}}}\Big[
\hat g(\xi\sin\theta)\hat f(\xi\cos\theta)-\hat g(0)\hat f(\xi)\Big]d\theta,
\end{align*}}
which provides \eqref{bob1}.
\end{proof}

\subsection{The harmonic oscillator}\label{6.sec.harmo}
The standard Hermite functions $(\phi_{n})_{n\geq 0}$ are defined for $x \in \rr$,
 \begin{multline}
 \phi_{n}(x)=\frac{(-1)^n}{\sqrt{2^n n!\sqrt{\pi}}} e^{\frac{x^2}{2}}\frac{d^n}{dx^n}(e^{-x^2})
 =\frac{1}{\sqrt{2^n n!\sqrt{\pi}}} \Bigl(x-\frac{d}{dx}\Bigr)^n(e^{-\frac{x^2}{2}})=\frac{ a_{+}^n \phi_{0}}{\sqrt{n!}},
\end{multline}
where $a_{+}$ is the creation operator
$$a_{+}=\frac{1}{\sqrt{2}}\Big(x-\frac{d}{dx}\Big).$$
The family $(\phi_{n})_{n\geq 0}$ is an orthonormal basis of $L^2(\R)$.
We set for $n\geq 0$, $\alpha=(\alpha_{j})_{1\le j\le d}\in\N^d$, $x\in \R$, $v\in \R^d,$
\begin{align}\label{}
\psi_n(x)&=2^{-1/4}\phi_n(2^{-1/2}x),\quad \psi_{n}=\frac{1}{\sqrt{n!}}\Bigl(\frac{x}2-\frac{d}{dx}\Bigr)^n\psi_{0},
\\
\Psi_{\alpha}(v)&=\prod_{j=1}^d\psi_{\alpha_j}(v_j),\quad \mathcal E_{k}=\text{Span}
\{\Psi_{\alpha}\}_{\alpha\in \N^d,\val \alpha=k},
\end{align}
with $\val \alpha=\alpha_{1}+\dots+\alpha_{d}$. The family $(\Psi_{\alpha})_{\alpha \in \nn^d}$ is an orthonormal basis of $L^2(\R^d)$
composed of the eigenfunctions of the $d$-dimensional harmonic oscillator
\begin{equation}\label{6.harmo}
\mathcal{H}=-\Delta_v+\frac{|v|^2}{4}=\sum_{k\ge 0}\Big(\frac d2+k\Big)\mathbb P_{k},\quad \text{Id}=\sum_{k \ge 0}\mathbb P_{k},
\end{equation}
where $\mathbb P_{k}$ is the orthogonal projection onto $\mathcal E_{k}$ whose dimension is $\binom{k+d-1}{d-1}$. The eigenvalue
$d/2$ is simple in all dimensions and $\mathcal E_{0}$ is generated by the function
$$\Psi_{0}(v)=\frac{1}{(2\pi)^{\frac{d}{4}}}e^{-\frac{\val v^2}{4}}=\mu_d^{1/2}(v),$$
where $\mu_d$ is the Maxwellian distribution defined in (\ref{maxwe}).

\subsection{On radial functions}\label{6.sec.radia}
If $u\in \mathscr S(\R^d)$ is a radial function
$$\forall x\in\rr^d,\,\,\forall A\in O(d),\quad  u(x)=u(Ax),$$
we define
$$\tilde u(t)=\frac{1}{\val{\SSS^{d-1}}}\int_{\SSS^{d-1}} u(t\sigma) d\sigma, \quad t\in \R.$$
This function is even, belongs to the Schwartz space $\mathscr S(\R)$ and satisfies
$$\forall t\in \R, \forall \sigma\in \SSS^{d-1}, \quad  \tilde{u}(t)=u(t\sigma), \quad \forall x\in \R^d, \quad  u(x)=\tilde{u}(\val x).$$
The Borel's theorem shows that the mapping $t\mapsto \tilde u(t)$ is also a Schwartz function of the variable $t^2$. We also recall that the Fourier transform of a radial function is radial and that the Fourier transformation is an isomorphism of the space of radial Schwartz functions.

\bigskip

\vs\noindent
{\bf Acknowledgements.}
The research of the second author was supported by the Grant-in-Aid for Scientific Research No.22540187, Japan Society for the Promotion of Science. The research of the third author was supported by the CNRS chair of excellence at Cergy-Pontoise University. The research of the last author was supported partially by ``The Fundamental Research Funds for Central Universities'' 
and the National Science Foundation of China No. 11171261.

\end{document}